\documentclass[11pt,leqno]{article}

\usepackage{amssymb}
\input psfig.sty

\newcommand{\maths}[1]{{\mathbb #1}}  

\newtheorem{defi}{Definition}[section]
\newtheorem{prop}[defi]{Proposition}
\newtheorem{theo}[defi]{Theorem}
\newtheorem{conj}[defi]{Conjecture}
\newtheorem{lemm}[defi]{Lemma}
\newtheorem{coro}[defi]{Corollary}
\newtheorem{rema}[defi]{Remark}
\newtheorem{exem}[defi]{Example}

\newcommand{\bdefi}{\begin{defi}}
\newcommand{\edefi}{\end{defi}}
\newcommand{\bprop}{\begin{prop}}
\newcommand{\eprop}{\end{prop}}
\newcommand{\btheo}{\begin{theo}}
\newcommand{\etheo}{\end{theo}}
\newcommand{\blemm}{\begin{lemm}}
\newcommand{\brema}{\begin{rema}}
\newcommand{\erema}{\end{rema}}
\newcommand{\bexer}{\begin{exem}}
\newcommand{\eexer}{\end{exem}}
\newcommand{\bconj}{\begin{conj}}
\newcommand{\econj}{\end{conj}}
\newcommand{\elemm}{\end{lemm}}
\newcommand{\bcoro}{\begin{coro}}
\newcommand{\ecoro}{\end{coro}}
\newcommand{\dem}{\noindent{\bf Proof. }}
\newcommand{\rem}{\noindent{\bf Remark. }}

\newcommand{\R}{{\cal R}}
\newcommand{\T}{{\cal T}}

\newcommand{\M}{{\cal M}}

\renewcommand{\H}{{\cal H}}

\newcommand{\ra}{\rightarrow}

\newcommand{\RR}{\maths{R}}
\newcommand{\NN}{\maths{N}}
\newcommand{\CC}{\maths{C}}
\newcommand{\QQ}{\maths{Q}}

\newcommand{\HH}{\maths{H}}

\newcommand{\ZZ}{\maths{Z}}

\newcounter{fig}

\def\myfigure#1#2#3{
\addtocounter{fig}{1}
\[
\begin{array}{c}
\mbox{\psfig{figure=#1.ps,height=#2}}\\
\\ 
\hbox{\rm Figure \arabic{fig} ~: #3.}
\end{array}
\]
}

\newcommand{\PSL}{\mbox{${\rm{ PSL}}_{2}(\CC)$}}
\newcommand{\SL}{\mbox{${\rm{ SL}}_{2}(\ZZ)$}}
\newcommand{\PSLZ}{\mbox{${\rm{ PSL}}_{2}(\ZZ)$}} 
\newcommand{\SLZI}{\mbox{${\rm{ SL}}_{2}({\ZZ}[i])$}}

\newcommand{\eop}[1]{{\hfill\fbox{\bf #1}}}

\title{Diophantine approximation\\ 
for
 negatively curved manifolds, I} 

\author{Sa'ar Hersonsky \and Fr\'ed\'eric Paulin}

\date{}

\begin{document}

\maketitle

\begin{abstract}
\noindent Let $M$ be a geometrically finite pinched negatively curved
Riemannian manifold with at least one cusp. Inspired by the theory of
diophantine approximation of a real (or complex) number by rational
ones, we develop a theory of approximation of geodesic lines starting
from a given cusp by ones returning to it. We define a new invariant
for $M$, the {\it Hurwitz constant} of $M$. It measures how well all
geodesic lines starting from the cusp are approximated by ones
returning to it. In the case of constant curvature, we express the
Hurwitz constant in terms of lengths of closed geodesics and their
depths outside the cusp neighborhood. Using the cut locus of the
cusp, we define an explicit approximation sequence for geodesic lines
starting from the cusp and explore its properties. We prove that the modular once-punctured hyperbolic torus has the minimum Hurwitz constant in its moduli space.\footnote{{\bf AMS codes:} 
53 C 22, 11 J 06, 30 F 40, 11 J 70. {\bf Keywords:} diophantine approximation, 
negative curvature, cusp, height of geodesics, Hurwitz constant.}
\end{abstract}

\section{Introduction}
\label{sec:Intro}

Let $M$ be a geometrically finite pinched negatively curved Riemannian
manifold with cusps. In this paper we study fine properties of the
geodesic flow of $M$ arising from the presence of cusps. There are
deep connections of diophantine approximation problems, for real or
complex numbers by elements of quadratic number fields, to hyperbolic
geometry (see \cite{Dan, For, For2, HV, HS, Pat, Ser3, Sch, Sul} and
the references therein).  The first purpose of this paper is to extend
the existing theory beyond the arithmetic case. We even allow the
curvature to be non constant.  To simplify the statements in the
intro\-duction, we assume that $M$ has finite volume and only one cusp
$e$.

A geodesic line starting from $e$ either converges to $e$, or
accumulates inside $M$. We say that the geodesic line is {\it
rational} in the first case, and {\it irrational} otherwise.  Let
$\HH^2$ denote the upper half space model of the hyperbolic plane.
Let $M$ be the orbifold $\HH^2/\PSLZ$.  In the rational (irrational) case,
the lift  starting from $\infty$ of the geodesic
line ends at a rational (irrational) point on
the real line.

For a rational line $r$, we introduce a new notion of complexity
(Definition~\ref{def:depth}), called the {\it depth} of $r$. The depth
$D(r)$ is the length of $r$ between the first and last meeting point
with the boundary of the maximal Margulis neighborhood of the
cusp. The set of depths of rational lines is a discrete subset of
$\RR$, whose asymptotic distribution will be exploited in a future
paper. There is a natural distance-like map $d$ (see section
\ref{subsec:dist_link}) on the space of geodesic lines starting from
the cusp. In constant curvature, it is given by the Euclidean metric
on the boundary of the maximal Margulis neighborhood of the cusp. In
general, it is a slight modification of the quotient of Hamenst\"adt's
metric \cite{Ham} on the horosphere covering that boundary.

The first goal of this paper is to give an analogue to the classical 
Dirichlet
theorem. When specialized to constant negative curvature, it coincides
with a known result for geometrically finite Kleinian
groups \cite{Pat, Ser1, HV}.

\btheo\label{theo:Const_intro} There
exists a positive constant $K$ such that for any irrational line $\xi$
starting from $e$, there exist infinitely many rational lines $r$ with
$$d(\xi,r) \leq K e^{-D(r)}.$$
\etheo

We call the infimum of such $K$ the {\it Hurwitz constant} of $e$, and
denote it by $K_{M,e}$. In the case that $M$ is the hyperbolic
orbifold $\HH^2/\PSLZ$, or the quotient of the hyperbolic $3$-space
$\HH^3$ by some {\it Bianchi group}, $K_{M,e}$ corresponds to the
classical Hurwitz constant for the approximation of real numbers by
rational ones, or the approximation of complex numbers by elements of
an imaginary quadratic number field.

The distribution (as a subset of $\RR$) of the Hurwitz constants of
one-cusped hyperbolic 3-manifolds or orbifolds is unknown, as is their
infimum. Among the one-cusped Bianchi 3-orbifolds, the smallest
Hurwitz constant is obtained by $\HH^3 /{\rm PSL}_2({\cal O}_{-3})$,
and afterwards seems to increase with the discriminant. 
\medskip

The second purpose of this paper is to give a geometric interpretation 
of the Hurwitz constant.  Using as height function a {\it Busemann
function} (giving the ``distance'' to the cusp, see
Section~\ref{sect:rational_rays}), normalized to have value $0$ on the
boundary of the maximal Margulis neigborhood (and to converge to
$+\infty$ in the cusp), we prove:

\btheo
The lower bound of the maximal
heights of the closed geodesics in $M$ is $-\log 2K_M$.
\etheo


In particular, there is a positive lower bound for the Hurwitz
constant of manifolds $M$ with a given upperbound on the length of the
shortest geodesic (since the length of a closed geodesic is at least
twice its height).

Let ${\cal M}_{g,1}$ be the moduli space of one-cusped hyperbolic
metrics on a connected oriented closed surface with genus $g\geq 1$
and one puncture. We prove (see Section \ref{sect:Hurw_torus}) that the
Hurwitz constant defines a proper continuous map from ${\cal M}_{g,1}$ to
$[0,+\infty[$. For higher genus, the one-cusped hyperbolic surfaces 
having the smallest Hurwitz constant in their moduli space are unknown.

\btheo\label{theo:Hurw_intro}
The Hurwitz constant is a proper real-analytic map on the orbifold ${\cal M}_{1,1}$, whose minimum $1/\sqrt{5}$ is attained precisely on 
the modular one-cusped hyperbolic torus.
\etheo

%

If $\Gamma$ is a geometrically finite Kleinian group which has
$\infty$ as a parabolic fixed point in the upper halfspace model, and
is normalized so that the boundary of the maximal Margulis
neighborhood of the cusp is covered by the horizontal plane through
$(0,0,1)$, one has the following result, where $c(\gamma)$ is the
lower left entry of the matrix of $\gamma$ in $\Gamma$. The formula
depends on the normalization.

\btheo
The Hurwitz constant of $\HH^3/\Gamma$ is
$$\inf_{\{\gamma\in\Gamma {\rm ~hyperbolic}\}}\; \frac{\sqrt{|{\rm
tr}^2 \gamma-4|}}{ \min_{\{\gamma'\in\Gamma\ {\rm is ~conjugate~to~}
\gamma\}} \;| 2c(\gamma')|}.$$ 
\etheo

Generalizing the definition in \cite{EP} for geometrically finite
hyperbolic manifolds with one cusp, we define the {\it cut locus of
the cusp} of $M$. This is the subset $\Sigma$ of points in $M$ from
which start at least two (globally) minimizing geodesic rays
converging to the cusp.  Let $L_e$ be the boundary of any standard Margulis
neighborhood of the cusp.  It follows that $M-\Sigma$ retracts in a
canonical way onto  $L_e$. Under some
technical assumptions (see section \ref{sect:cut_locus}) which are
satisfied in the constant curvature case, and that will be assumed
throughout the introduction, $\Sigma$ has a natural smooth locally
finite stratification.
\medskip

The third purpose of this paper is to construct an explicit good
approximation to any irrational line $\xi$ starting from $e$ by
rational ones.  Using the cut locus of the cusp $\Sigma$, we define in
section \ref{sect_good_approximation} an explicit sequence of rational
lines $r_n$ converging to $\xi$. To simplify the definition, assume
here that $\xi$ is transverse to the stratification of $\Sigma$.  In
particular, $\xi$ cuts $\Sigma$ in a sequence of points $(x_n)_{n\in
\NN}$ with $x_n$ belonging to an open top-dimensional cell $\sigma_n$
of $\Sigma$.  For each $n\in \NN$, let $c_n$ be the path consisting of
the subsegment of $\xi$ between $e$ and $x_n$, followed by the
minimizing geodesic ray from $x_n$ to the cusp $e$, starting on the
other side of $\sigma_n$.  Let $r_n$ be the unique geodesic line
starting from $e$, properly homotopic to $c_n$.
\medskip

\btheo The geodesic line
$r_n$ is rational and there exists a constant $c>0$ (independant of
$\xi$) such that for every $n$ in $\NN$,
$$d(r_n,\xi)\leq c \;e^{-D(r_n)}.$$ 
\etheo

For instance, if $\gamma$ is a closed geodesic transverse to the
stratification of $\Sigma$, then the geodesic line $\xi$ starting from
the cusp and spiraling around $\gamma$ has an eventually periodic
sequence $(r_n)_{n\in\NN}$. Conversely, if a geodesic line $\xi$,
starting from the cusp and transverse to the stratification of
$\Sigma$, has an eventually periodic sequence $(r_n)_{n\in\NN}$, then
$\xi$ spirals around a closed geodesic in $M$.

\medskip
The fourth purpose of this paper is to express an endpoint of an
irrational line by means of geometrical data obtained from the
manifold. The cut locus of the cusp $\Sigma$ allows us to parametrize
the geodesic lines starting from the cusp by sequences
$(a_n)_{n\in\NN}$, where the $a_n$'s belong to a countable alphabet
(see section \ref{sect:continued_fraction}).  This alphabet is the set
$\pi_1(L_e,\R)$ of homotopy classes of paths in $L_e$ (relative
boundary) between points of a geometrically defined finite subset
$\R$ of $L_e$. An irrational line $\xi$ which is transverse to
$\Sigma$ travels from one point $x_n$ of $\Sigma$ to its next
intersection point $x_{n+1}$ with $\Sigma$. By the property of
$\Sigma$, the subpath of $\xi$ between $x_n$ and $x_{n+1}$ is
homotopic to a path lying on $L_e$, which is $a_{n+1}$ (see section
\ref{sect:continued_fraction}).

We prove (Theorem \ref{theo:continued_fraction_determines_ray})
that $\xi$ is uniquely determined by the sequence
$(a_n)_{n\in\NN}$. In constant curvature and dimension $3$, where
$L_e$ is a $2$-torus, we give an explicit formula (Theorem
\ref{theo:new_formula}) giving the $r_n$'s in terms of the $a_n$'s,
analogous to the expression giving the $n$-th convergent of an
irrational real number in terms of its continued fraction
developement.  Our search for an explicit formula expressing the
endpoint of an irrational line by geometrical data was inspired by
\cite[Theorem A]{Ser3} (the \SL\ case) and \cite[Part II]{For} (the
\SLZI\ case).  But, even in the special cases above, our formula gives
new information and calls for further study, among others of the
growth of the depths of rational lines in the good approximation
sequence of an irrational line. One should also consult the work in
\cite[chapters 2-5]{Sch}, giving a continued fraction expansion for
some complex numbers. Schmidt's construction is completely different
from ours.

\medskip
Two more papers in this series are under preparation.  In the
second we will give a coding of the geodesic flow in $M$, using
cutting sequences of all geodesic lines with the dual tessalation of
the cut locus of the cusp, by a subshift of finite type on a countable
alphabet.  In the third we give an analogue of the Khinchine-Sullivan
theorem, and an estimate of the Hausdorff dimension, in terms of $s$
and the bounds on the curvature, of the set of geodesic lines $\xi$
starting from a cusp $e$ for which there exists infinitely many
rational lines $r$ with $d(r,\xi)\leq e^{-sD(r)}$.

\bigskip
\noindent {\small {\bf Acknowledgement:} We are indebted to Steven
Kerckhoff, Darren Long, Nikolai Makarov and Alan Reid for their
encouragment, support and many valuable conversations. We are grateful
to Dave Gabai for his encouragment.  We thank Peter Shalen for telling us
about Vulakh's work.  Parts of this paper were
written during a stay of the first author at the IHES in August-September
1998, and during a one week stay of the second author at the Warwick
University under an Alliance project. We thank Caroline Series for
many helpful comments.}

\section{Rational and irrational lines starting from a cusp}
\label{sect:rational_rays}

Let $M$ be a (smooth) complete Riemannian $n$-manifold with pinched
negative sectional curvature $-b^2\leq K\leq -a^2<0$.  Fix a universal
cover $\widetilde{M}$ of $M$, with covering group $\Gamma$.

A {\it geodesic segment, ray, line} in a $M$ is a locally isometric
map from a compact interval, a minorated unbounded interval, $\RR$
respectively, into $M$.  Note that any geodesic segment, ray, line in
$\widetilde{M}$ is (globally) minimizing, but that it is not always
the case in $M$.

The boundary $\partial \widetilde{M}$ of $\widetilde{M}$ is the space
of asymptotic classes of geodesic rays in $\widetilde{M}$. Endowed
with the cone topology, the space $\widetilde{M}\cup
\partial\widetilde{M}$ is homeomorphic to the closed unit ball in
$\RR^n$ (see for instance \cite[Section 3.2]{BGS}). The {\it limit
set} $\Lambda(\Gamma)$ is the set $\overline{\Gamma x}\cap
\partial\widetilde{M}$, for any $x$ in $\widetilde{M}$.  See  
\cite{Bow} for the following definitions.


\bdefi With $\Gamma$ and $M$ as above:
\begin{enumerate} 
\item A point $\xi$ in $\Lambda(\Gamma)$ is a {\rm conical limit point} of
$\Gamma$ if it is the endpoint of a geodesic ray in $\widetilde{M}$
which projects to a geodesic in $M$ that is recurrent in some compact
subset.
\item 
A point $\xi$ in $\Lambda(\Gamma)$ is a {\rm bounded parabolic point}
if it is fixed by some parabolic element in $\Gamma$, and if the
quotient $(\Lambda(\Gamma)-\{\xi\})/\Gamma_{\xi}$ is compact, where
$\Gamma_{\xi}$ is the stabilizer of $\xi$.
\item 
The group $\Gamma$ and the manifold $M$ are called {\rm geometrically
finite} if every limit point of $\Gamma$ is conical or bounded
parabolic.
\end{enumerate}
\edefi

We assume in this paper that $M$ is geometrically finite and {\it non
elementary}, i.e.~that the limit set contains at least $3$ (hence
uncountably many) points. In that case, the limit set is the smallest
non empty invariant closed subset of $\partial \widetilde{M}$.  The
{\it convex core} $C(M)$ of $M$ is the image by the covering map
$\widetilde{M}\ra M$ of the convex hull of the limit set of
$\Gamma$. For instance, if $M$ has finite volume, then $M$ is
geometrically finite and $C(M)=M$.

A {\it cusp} in $M$ is an asymptotic class of  minimizing
geodesic rays in $M$ along which the injectivity radius goes to $0$.
If $M$ has finite volume, the cusps are in one-one corespondance with 
the ends of $M$. 
A geodesic ray (line) {\it converges to} the cusp if some positive
subray is asymptotic to a ray in the equivalence class of the cusp. A
geodesic ray converges to some cusp if and only if some (any) lift in
$\widetilde{M}$ ends in a parabolic fixed point.  In all that follows,
we fix a cusp $e$.

Given any minimizing geodesic ray $r$, recall (see \cite[Section
3.3]{BGS}) that the Buseman function $\beta_r:M\ra \RR$ is the
$1$-Lipschitz map defined by the limit (which exists for all $x\in M$)
$$\beta_r(x)=\lim_{t\ra\infty} \bigl(t-d_M(x,r(t))\bigr).$$ 

Fix a minimizing ray $r$ converging to the cusp $e$. Let
$\widetilde{r}$ be any lift of $r$ to $\widetilde{M}$. The following
facts follow from the Margulis Lemma (see for instance \cite{BK} and
\cite[Sections 9-10]{BGS}) and since $M$ is non elementary. There
exists $\eta_0=\eta_0(r)$ in $\RR$ such that, given $t$ in $\RR$, 
the quotient of
the horosphere $\beta_{\widetilde{r}}^{-1}(t)$ by the stabilizer in
$\Gamma$ of the endpoint $\widetilde{r}(+\infty)$ of $\widetilde{r}$
embeds in $M$ under the covering map $\widetilde{M}\ra M$ if and only
if $t>\eta_0$.  So that for $\eta>\eta'> \eta_0$, the level set
$\beta_r^{-1}(\eta)$ identifies with the quotient of an horosphere
$\beta_{\widetilde{r}}^{-1}(\eta)$ by the stabilizer
$\Gamma_{\widetilde{r}(+\infty)}$. There is a unique minimizing
geodesic ray starting perpendicularly to $\beta_r^{-1}(\eta')$ at a
given point $x\in\beta_r^{-1}(\eta')$, and entering
$\beta_r^{-1}(]\eta',+\infty[)$. It converges to $e$ and meets
$\beta_r^{-1}(\eta)$ perpendicularly in exactly one point
$\phi(x)$. The distance $d_M(x,\phi(x))$ is the constant
$\eta-\eta'$. If the curvature is constant $-1$, then the
homeomorphism $\phi:\beta_r^{-1}(\eta') \ra \beta_r^{-1}(\eta)$
induces a contraction of the induced length metrics of ratio
$e^{\eta'-\eta}< 1$. 

Define $\beta_{e}(x) = \beta_r(x)-\eta_0$. Since Busemann functions of
asymptotic minimizing rays differ by an additive constant, the map 
$\beta_{e}:M\ra\RR$ does not depend on $r$.

\bdefi({\bf Busemann function of the cusp}). 
\label{def:normalizedBu}
The map $\beta_{e}:M\ra\RR$ is called the Busemann function
of $e$, and $\beta_e^{-1}(]0,+\infty[)$ the maximal Margulis
neighborhood of $e$
\edefi


\medskip
Consider the set of geodesic lines $c:\,]-\infty,+\infty[\,\ra M$
whose negative subrays converge to $e$ and that are recurrent in some
compact subset or ending in some cusp. Identify two of them if they
differ by a translation of the time. An equivalence class will be
called a {\it geodesic line starting from $e$}. The set of equivalence
classes will be called the {\it link of $e$ in $C(M)$}, and denoted by
$Lk(e,C(M))$. It is in one-to-one correspondance with a closed subset 
of the intersection
$\beta_e^{-1}(\eta)\cap C(M)$ of the convex core and a level set of
the Busemann function of $e$, for any $\eta>0$, by the map which
associates to $c$ its first intersection point with
$\beta_e^{-1}(\eta)$.

\bdefi ({\bf Rational and irrational geodesic lines}).  A geodesic
line starting from $e$ will be called {\rm rational} if it converges
to $e$. A geodesic line starting from $e$ which is not rational and
which does not converge into another cusp  will be called {\rm irrational}.
\edefi

Being a rational line is equivalent to requiring that the line meets
perpendicularly a second time a level set $\beta_e^{-1}(\eta)$, for
any $\eta>0$ (it is contained after that time in
$\beta_e^{-1}([\eta,+\infty[)$ and converges to $e$).

\subsection{The distance on the link of the cusp}
\label{subsec:dist_link}

For our approximation purpose, we will need to measure how close are
two points in the link of $e$ in $C(M)$. For that we will define a
``distance-like'' map $d:Lk(e,C(M))\times Lk(e,C(M))\ra [0+\infty[$ as
follows.

Let $H$ be any horosphere in $\widetilde{M}$ centered at
$a\in\partial\widetilde{M}$. We first define a map $d_H:H\times H\ra
[0+\infty[$. Let $x$ be a point on $H$, let $L_x$ be the geodesic line
through $x$ starting at $a$ and let $x'$ be its endpoint (see Figure
\addtocounter{fig}{1}\arabic{fig}\addtocounter{fig}{-1}), so that
$L_x$ is oriented from $x$ towards $x'$. For $r>0$, let $H_r$ be the
horosphere centered at the endpoint of $L_x$, meeting $L_x$ at a point
$u$ at signed distance $-\log 2r$ of $x$ along $L_x$.  For every $x,y$
in $H$, define $d_H(x,y)$ to be the infimum of all $r>0$ such that
$H_r$ meets $L_y$.

\myfigure{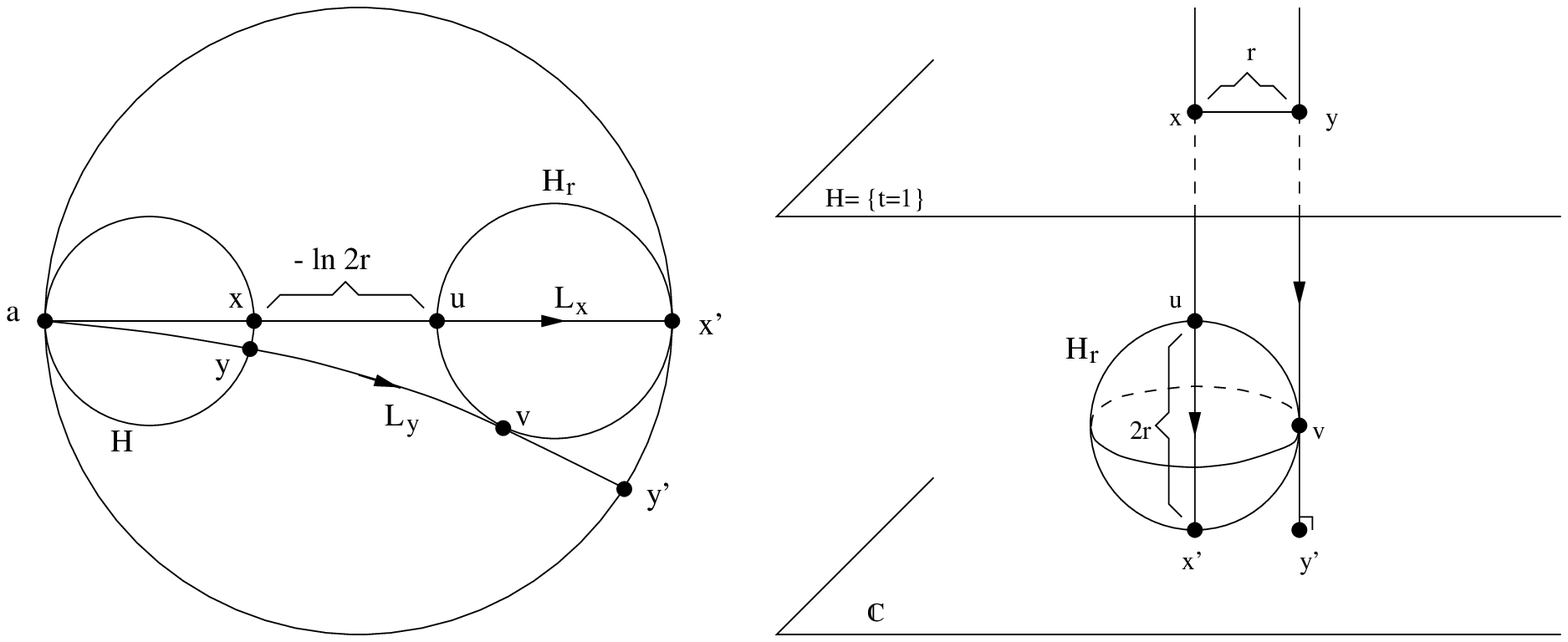}{5.8cm}{The ``metric'' on the link of a cusp}

This map $d_H$ is a priori not symmetric nor transitive (see
\cite[Appendix]{HP} for a related actual distance). In constant
curvature, it coincides with the induced Riemannian metric (which is
flat) on the horosphere $H$. To prove that (see Figure \arabic{fig}),
by naturality of the construction, one may assume that $a$ is the
point at infinity in the upper half-space model with curvature $-1$,
and $H$ is the horizontal horosphere $t=1$. Let $r$ be the Euclidean
distance between $x$ and $y$. Let $H_r$ be the horosphere which is
tangent to the horizontal coordinate hyperplane at the vertical
projection $x'$ of $u$, and has Euclidean radius $r$. It bounds a
horoball which is the smallest one meeting the vertical line through
$y$. An easy computation shows that the hyperbolic distance between
$x$ and $u$ is $-\log 2r$.

If for any two points on an horosphere $H$ in $\widetilde{M}$, there
exists an isometry in $\widetilde{M}$ preserving the horoball bounded
by $H$ which exchanges the two points, then $d_H$ is symmetric. This
is the case if $\widetilde{M}$ is a symmetric space (of non compact
type) of rank one.

From the topological point of view, the map $d_H$ is as good as a
distance:

\bprop\label{prop:distance_link} 
For every $x$ in $H$, let $B_{d_H}(x,\epsilon)=\{y\in H\,|\,d_H(x,y)
<\epsilon\}$. Then $\{B_{d_H}(x,\epsilon)\,|\,\epsilon>0\}$ is a
fundamental system of neighborhoods at $x$.  
\eprop

\dem Since $\widetilde{M}$ is a negatively curved Riemannian manifold,
the map $H\ra \partial\widetilde{M}-\{a\}$ which sends $x$ to $x'$ is
a homeomorphism.  So we only have to prove that
$\{B'(x',\epsilon)\,|\,\epsilon>0\}$ is a fundamental system of
neighborhoods at $x'$, with $B'(x',\epsilon)=\{y'\in
\partial\widetilde{M}-\{a\}\,|\,d_H(x,y)<\epsilon\}$.  Assume that
$y'\neq x'$, let $HB'$ be the smallest horoball centered at $x'$ and
meeting $L_y$. Let $H'$ be the boundary of $HB'$, then $H'$ is tangent
at $L_y$ in a point $v$ and meets $L_x$ in a point $u$ (see Figure
\arabic{fig}). Let $p$ be the perpendicular projection of $v$ on
$L_x$. By convexity of the horoballs, the point $p$ belongs to $HB'$,
so that $x,u,p$ are in this order on $L_x$. By the existence of a
negative upper bound on the curvature, there exists (see for instance
\cite{GH}) a constant $C>0$ (depending only on the bound) and a map
from $L_x\cup L_y$ into a tree $T$ which is an isometry on $L_x$ and
on $L_y$, and preserves the distances up to the additive constant $C$.
Let $\overline{z}$ be the image in $T$ of any $z\in L_x\cup L_y$, and
$\overline{a},\overline{x'}, \overline{y'}$ the ends of $T$
corresponding to $a, x',y'$. Take $\overline{x}$ as basepoint in $T$.
Since $L_x,L_y$ are asymptotic in $a$, their images
$\overline{L_x},\overline{L_y}$ meet in a subray from $P$ to the end
$\overline{a}$, where $P$ is the point in $T$ which is the projection
of $\overline{y'}$ on $\overline{L_x}$.  The map $z\mapsto
\overline{z}$ preserves the distance (and hence also the Busemann
functions) up to the additive constant $C$. By the properties of the
Busemann functions in trees, since $u,v$ lie on the same horosphere in
$\widetilde{M}$, the points $\bar{u},\bar{v}$ are on horospheres centered 
at $\bar{x}$ at distance at most $2C$. So that  if $\overline{v}$ 
lies between $P$ and $\overline{y'}$, then
$$|d_T(\overline{u},P)-d(\overline{v},P)|\leq 2C$$
and if $\overline{v}$ does not lie between $P$ and $\overline{y'}$, 
then both $\overline{u}$ and $\overline{v}$ lie on $\overline{L_x}$, 
hence 
$$d_T(\overline{u},\overline{v})\leq 2C.$$ Since $u$ lies on the smallest horosphere centered at $x'$ and meeting $L_y$, and since the smallest horoball
in $T$ centered at $\overline{x'}$ which contains a point of
$\overline{L_y}$ has its horosphere passing through $P$, it follows that
$d_T(\overline{u},P)\leq 2C$. Hence if $\overline{v}$ lies
between $P$ and $\overline{y'}$, one has $$d_T(\overline{u},
\overline{v})=d_T(\overline{u},P)+d(\overline{v},P)\leq 6C.$$ 
Therefore the points $u$ and $v$ are at distance in
$\widetilde{M}$ at most $7 C$.  Hence $x'$ and $y'$ are close on
$\partial\widetilde{M}$ if and only if $\overline{x'}, \overline{y'}$
are close on $\partial T$, which is equivalent to $\overline{x},
P,\overline{x'}$ being in this order on $\overline{L_x}$ and
$d_T(\overline{x}, P)$ is big, which occurs if and only if $x, u,x'$ are 
in this
order on $L_x$ and $d_{\widetilde{M}}(x, u)$ is big. This proves the
result.  \eop{\ref{prop:distance_link}}

\bigskip
We now define a map $d$ on Lk$(e,C(M))$ by taking quotients.

\bdefi
\label{def:dist_on_the_link}
For any $\eta>0$, let $H_\eta$ be any horosphere in $\widetilde{M}$
covering $\beta_e^{-1}(\eta)$, and let $\pi:\widetilde{M}\ra M$ be the
covering map. For any two points $x,y$ in Lk$(e,C(M))$, identified as
above with a closed subset of $\beta^{-1}(\eta)$, define $d(x,y)$ as the 
infinum
of $e^{\eta}d_{H_\eta}(\tilde{x},\tilde{y})$ for all the preimages
$\tilde{x}, \tilde{y}$ of $x,y$ respectively.  
\edefi

It is easy to show that the map $d$ on Lk$(e,C(M))\times$ Lk$(e,C(M))$
does not depend on the choice of the horosphere $H_\eta$ (by
equivariance), nor on $\eta>0$. This is because if $H,H'$ are
horospheres in $\widetilde{M}$ centered at the same point at infinity,
with say $H'$ contained in the horoball bounded by $H$, with $\Delta$
the (constant) distance between $H'$ and $H$, then
$d_{H'}=e^{-\Delta}d_H$.

In general, $d$ is not a distance. It is a distance if $M$ is {\it
locally symmetric}. If furthermore the curvature is constant $-1$,
then $d$ coincides with the (flat) induced length metric on
$\beta^{-1}(\eta)$ normalized by $e^{\eta}$.

\medskip
For $X$ a complete simply connected Riemannian manifold with
pinched negative curvature with the upper bound on the curvature being
$-1$, Hamenstadt defined (see \cite{Ham}) a metric on the space at
infinity minus a point $a$, depending on an horosphere $H$ centered at
$a$, as follows. The distance between $x,y$ in $\partial X-\{a\}$ is given by
$$\delta(b,c)=\lim_{t\ra \infty} e^{-t +\frac{1}{2}d_X(b_t,c_t)}$$ where
$t\mapsto b_t, t\mapsto c_t$ are geodesic lines starting from $a$,
ending in $b,c$ and passing through $H$ at time $0$. In
\cite[Appendix]{HP}, we proved that this limit exists and 
defines a metric for any CAT$(-1)$ space $X$.

By identifying $\partial X-\{a\}$ with the horosphere $H$ as usual
(sending $b\neq a$ to the point of intersection with $H$ of the
geodesic line starting at $a$ and ending at $b$), one gets a metric
$\delta_H$ on $H$. Our metric $d_H$ is equivalent to $\delta_H$:

\brema 
Assume that the upper bound of the curvature of $M$ is $-1$.
Then there exists a constant $c>0$ such that for every $x,y$ in $H$
$$\frac{1}{2}\delta_H(x,y)\leq d_H(x,y)\leq c \delta_H(x,y).$$ 
\erema

\dem 
Keeping the notations of the definition of $d_H$ (and Figure 1),
we have, since $u$ (resp.~$v$) lies on the segment between $x$
(resp.~$y$) and $x_t$ (resp.~$y_t$) for $t$ big enough, by the triangle
inequality, since the minimum of the distances between a point of $H$
and a point of $H_r$ is attained by $d_X(x,u)$, and since $u,v$ lie on
the same horosphere centered at the point to which $x_t$ converges as
$t\ra\infty$, if $\epsilon$ is small enough, then for $t$ big enough,
$$\begin{array}{cl}
-2t +d_X(x_t, y_t)  & =
- [d_X(x, u) + d_X(u, x_t)] - [d_X(y,v) + d_X (v,y_t)] + d_X(x_t, y_t) 
\\ &
\leq 
- d_X(x, u) - d_X(u, x_t) - d_X(y,v) - d_X (v,y_t)
+ d_X(x_t, v) + d_X(v,y_t)
\\ &
\leq -2 d_X(x, u) + [d_X(v,x_t) -d_X(x_t,u)]
\\ &
\leq -2 d_X(x, u) + \epsilon
= 2 \log ( 2 d_H(x,y)) + \epsilon
\end{array}$$
so that 
$$\delta_H(x,y)\leq 2 d_H(x,y).$$ 
By the technique of approximation by trees, one may also show that there 
is an explicitable universal constant $c$ such that
$$d_H(x,y)\leq c \delta_H(x,y).$$
\eop{}

\medskip
Since $\delta_H$ is a distance (inducing the right topology), this
gives another proof of proposition \ref{prop:distance_link}.

\subsection{The depth of rational geodesic lines}
\label{subsect:depth}

In this subsection we define a notion of complexity for rational
geodesic lines. First we present a connection between rational
geodesic lines, which were defined in geometric terms and a set of
double cosets in $\pi_1(M)$. This will allow us to perform
computations in the constant negative curvature case.

\medskip
Choose a base point on the level set $L=\beta_e^{-1}(1)$. (The first
intersection point of the geodesic line starting from $e$ passing
through the base point with any level set $\beta_e^{-1}(\eta)$ for
$\eta >0$ gives a base point on that level set and we will use the
subsegments of that geodesic line to identify the fundamental groups
of $M$ based at the base points on different level sets.)  By the
Margulis Lemma, the inclusion $i:L\ra M$ induces an injection between
fundamental groups $i_{\star}:\pi_1 L\ra\pi_1 M$ (use the choosen base 
point both for $L$ and $M$). We identify $\pi_1 L$ with its image.

\blemm
\label{lem:identify}
The set of rational lines is in one-to-one correspondence with the set
of double cosets $\pi_1 L\backslash\pi_1 M/\pi_1 L$.  
\elemm

\dem 
If $r$ is a rational line, and $c$ is its subpath between the two
succesive perpendicular intersection points $x,y$ with $L$, then
choosing a path on $L$ between the base point and $x$, and from $y$ to
the base point defines an element in $\pi_1 M$ whose double coset is
uniquely defined. Conversely, the straightening process (inside a
given homotopy class of paths in $M- \beta_e^{-1}(]1,+\infty[)$ with
endpoint staying in $L$) associates a rational line to each double
coset. More precisely, to any path $c$ in $M-
\beta_e^{-1}(]1,+\infty[)$ with endpoints on $L$, let $\widetilde{c}$
be a lift of $c$ to $\widetilde{M}$.  Its endpoints belongs to two
lifts $\widetilde{L}_1, \widetilde{L}_2$ which are disjoint, unless
$c$ is homotopic into $L$. The horoballs bounded by $\widetilde{L}_1,
\widetilde{L}_2$ are closed convex subsets of $\widetilde{M}$ with no
common point at infinity. Since the curvature is non positive, there
exists a unique common perpendicular segment $s$. The projection of
$s$ to $M$ gives a geodesic segment in $M$ homotopic to $c$ by an
homotopy moving the endpoints along $L$, which is a subsegment of a
rational line.  Clearly, the above two maps are inverse
one of the other.  \eop{\ref{lem:identify}}

\medskip
We measure the complexity of a rational line $r$ by the length
$\ell_\eta(r)$ of its subsegment between the two perpendicular
intersection points with any level set $\beta_e^{-1}(\eta)$ for
$\eta>0$, suitably normalized not to depend on $\eta$.

\bdefi ({\bf Depth of a rational line}).
\label{def:depth}  
The depth of a rational line $r$ is $D(r)= \ell_\eta(r) -2\eta$.
\edefi

By the properties of the level sets $\beta_e^{-1}(\eta)$, the depth
$D(r)$ is independant of $\eta$.

\brema \label{rem:non_numerote}
The set of depths of rational lines starting from $e$ is a
discrete subset of $\RR$, with finite multiplicities.  
\erema

Indeed, the set of intersections with the preimage of $C(M)$ of the horospheres
covering $\beta_e^{-1}(1)$ is locally finite in the universal covering
of $M$, since the group acts discretely, and the stabilizer of each
such intersection acts cocompactly on it.

\subsection{The constant curvature case}
\label{sect:const_Curv_Case}

In the case of constant negative curvature, we have the following
precise description, which also extends to the orbifold case. We
restrict to the dimension $3$, though everything is valid in higher
dimensions using the Vahlen matrices (see for instance \cite{Ah}).

Let $M$ be a connected orientable geometrically finite complete 
hyperbolic $3$-orbifold. Then, according to the description of
the thin part of $M$, each cusp $e$ of $M$ has a neighborhood $N$
isometric (for the induced length metric) to $(\T\times [a,+\infty[,ds^2)$ 
where $a=a(N)\in \RR$ is some constant, $\T=\T(e)$ is a
connected orientable flat $2$-orbifold with metric $dx^2$,
and $ds^2=e^{-2t}dx^2 + dt^2$. In all what follows, the cusp $e$ is fixed.

We will call such an $N$ a {\it standard cusp neighborhood} (of $e$).
The union of all standard cusps neighborhood is isometric to ${\cal
T}\times ]a_0,+\infty[$ endowed with the metric $ds^2=e^{-2t}dx^2 +
dt^2$ for some $a_0>0$, and we fix such an isometry.  Choose a base
point $\ast$ on ${\cal T}$. This gives a choice of base-point as above
by considering the geodesic ray corresponding to $\{\ast\}\times
]a_0,+\infty[$. Using orbifold fundamental groups, the set of rational
lines is in one-to-one correspondence with the set of double cosets
$\pi_1^{{\rm orb}}\partial N\backslash\pi_1^{{\rm orb}} M/\pi_1^{{\rm
orb}}\partial N$.

We note that for a geodesic $c$ in an orbifold that meets the singular
locus of the orbifold $M$ (which is a finite metric graph with
possibly some points removed, corresponding to some ends) at a point
$x$, then ingoing and outgoing tangent vectors of $c$ at $x$ make an
angle strictly less than $\pi$ (unless $c$ is locally contained in the
singular locus).

\medskip
Using a suitable uniformization of $M$, the depth of a rational line
can be computed in algebraic terms. We will use the upper half-space
model $\{(z,t)\,|\, z\in \CC, t>0\}$ for the real hyperbolic $3$-space
$\HH^3$, with the metric $$ds^2= \frac{|dz|^2 +dt^2}{t^2}.$$ We fix an
isometry between $M$ and $\HH^{3}/\Gamma$ with $\Gamma$ a subgroup of
Isom$_{+}(\HH^{3})$, in the following way.

By the Cartan-Hadamard theorem, there exists such an isometry, which
induces an isomorphism between the orbifold fundamental group
$\pi_1^{{\rm orb}} M$ of $M$ and $\Gamma$. Up to conjugating $\Gamma$,
we may assume that the stabilizer of the point at infinity in the
upper half space model, $\Gamma_\infty$ corresponds to $\pi_1^{{\rm
orb}} (\partial N)$, for $N$ any standard cusp neighborhood of
$e$. For any $h>0$, let $H_\infty(h)$ be the {\it horizontal
horosphere}, defined by the equation $t=h$ in the upper
half-space. Note that $\Gamma_\infty$ is a discrete group of
isometries of each $H_\infty(h)$ in $\HH^3$.

We may assume, up to conjugating $\Gamma$ by a dilatation, that
$H_\infty(h)/\Gamma_\infty$ injects into $M$ if and only if $h>1$
(this is coherent with the previous normalisation).  So that each
standard cusp neighborhood of $e$ is of the form
$H^+_\infty(h)/\Gamma_\infty$ for some $h>1$, where $H^+_\infty(h)$ is
the set of points $(z,t)$ with $t>h$. The uniformization $\HH^3\ra M$
is well defined only up to precomposition by an isometry of the form
$z\mapsto a z+b$ with $a,b$ complex and $|a|=1$, but the subsequent 
constructions will not depend on the choice of the uniformization.

In all what follows, we will identify the link of $e$ in $M$ with
$\CC/\Gamma_\infty $, by the map which sends a geodesic line $r$
starting from $e$ to the endpoint on the horizontal coordinate plane
of any lift $\widetilde{r}$ starting at $\infty$ of $r$, modulo
$\Gamma_\infty$. The distance $d$ on $Lk(e,M)$ defined in
Definition~\ref{def:dist_on_the_link} is exactly the Euclidean
distance on $\CC/\Gamma_\infty$, since the map $(z,1)\mapsto z$ from
$H_\infty(1)$ (endowed with the induced length metric) to $\CC$ (with
the Euclidean metric) is an isometry.

Since Isom$_+(\HH^3) = \PSL =$ SL$_2(\CC)/\{\pm 1\}$, we will write each
$\gamma\in\Gamma$ as
$$\gamma=\pm\left(\begin{array}{cc} a(\gamma) & b(\gamma)\\ c(\gamma)&
d(\gamma)\end{array}\right),$$ which acts on $\CC\cup\{\infty\}$ as
usual by the M\"obius transformation $\displaystyle{z\mapsto
\frac{a(\gamma)z+b(\gamma)}{c(\gamma)z+d(\gamma)}}$.

\medskip
Using these identifications we obtain an explicit expression for the
endpoint (on the complex plane) of any rational ray $r$.

\blemm\label{lem:depth_computation}
Let $\Gamma_\infty\gamma\Gamma_\infty$ be the double coset
associated to a rational line $r$. Then $r = \frac{a(\gamma)}{c(\gamma)} 
{\rm ~mod~} \Gamma_\infty$ and $D(r)=2\log |c(\gamma)|$.
\elemm

\dem 
Let $\widetilde{r}$ be any lift of $r$ to $\HH^3$ starting from
$\infty$, which is a descending vertical geodesic. Any two such lifts
differ by the action on the left of an element of $\Gamma_\infty$. By
the definition of $\gamma$, there exists an element $\alpha$ in
$\Gamma_\infty$ such that the endpoint of $\widetilde{r}$ on the
horizontal coordinate plane is $\alpha\gamma(\infty)$.  Since $\gamma$
does not fix $\infty$, $c(\gamma)$ is different from zero. For any
$\alpha$ in $\Gamma_\infty$ we have $\alpha\gamma(\infty) =
\alpha(\frac{a(\gamma)}{c(\gamma)})$. By taking the image in
$\CC/\Gamma_\infty$ the first assertion of the lemma follows.

The horizontal horosphere $H_\infty(1)$ centered at $\infty$ is mapped
by $\gamma$ to the horosphere centered at $\gamma(\infty)$ of
Euclidean diameter $\frac{1}{|c(\gamma)|^2}$. By the definition of the
depth, this implies the second assertion of the lemma (which is
independant of the choice of the representative of the double coset).
\eop{\ref{lem:depth_computation}}

\section{The approximation constant}
\label{sect:Hurwitz}

In this section we study the approximation of irrational rays by
rational ones, in pinched negatively curved Riemannian manifolds. We
keep the same notation as in the beginning of the previous section, in
particular for the distance $d$ on $Lk(e,C(M))$.

Theorem~\ref{theo:existence_Hurwitz} establishes an analogue of the
classical approximation theorem by
Dirichlet. Theorem~\ref{theo:positivity_Hurwitz} gives an analogue of
the classical Hurwitz constant. These two results were known in
constant curvature (see for example
\cite{HS,HV,Vul1,Vul2}). Theorem~\ref{theo:equality_four_terms} gives
a relation between the Hurwitz constant and the {\it heights} of
closed and non-cusp converging geodesic lines, which is new even for
Fuschian or Kleinian groups. Corollary~\ref{coro:Hurwitz_constcurv}
expresses the Hurwitz constant for a Kleinian group $\Gamma$ 
in algebraic terms, once a natural
normalization of $\Gamma$ has been made. At the end of this section we
present some known values of the Hurwitz constant of the Bianchi
groups.

\btheo \label{theo:existence_Hurwitz} 
Let $M$ be a pinched negatively curved geometrically finite Riemannian
manifold. For any irrational line $\xi$ starting from a cusp $e$,
there exists a constant $K>0$ such that for infinitely many rational
lines $r$, one has
$$ d(\xi,r) \leq K e^{-D(r)}.$$
\etheo

For a given $\xi$, the infimum over all such $K$'s is denoted by
$K(\xi)$, and is called the {\it Hurwitz constant} of $\xi$. The
supremum of $K(\xi)$ over all irrational lines $\xi$ will be called
the {\it Hurwitz constant} of $e$, and will be denoted by $K_{M,e}$.

\btheo \label{theo:positivity_Hurwitz}
Let $M$ be a pinched negatively curved geometrically finite Riemannian
manifold. For every cusp $e$ of $M$, we have
$$0<K_{M,e}<\infty.$$ 
\etheo

The maximum of $K_{M,e}$ on the finitely many cusps $e$ of $M$ is a
new invariant of the geometrically finite pinched negatively curved
manifold $M$, that we call the {\it Hurwitz constant of $M$}. 
Theorem \ref{theo:Const_intro} of the introduction follows from 
Theorems \ref{theo:existence_Hurwitz} and \ref{theo:positivity_Hurwitz}

Let $P$ be a non empty closed subset of $C(M)$ which does not meet
some neighborhood of $e$. Define the {\it height of $P$ with respect
to the cusp $e$} to be the maximum of the normalized Busemann function
of $e$ on $P$, that is
$$ht(P)=\max_{x\in P} \;{\beta_{e}(x)}.$$ 

To prove that the maximum is attained, set $t=\sup_{x\in
P}\beta_{e}(x)\in ]0,+\infty[$. Since $M$ is geometrically finite and
$P$ is contained in $M$, the subset
$P\cap(\beta_e^{-1}([t-1,+\infty[)-\beta_e^{-1}(]t+1,+\infty[))$ is
compact. The maximum of $\beta_{e}$ on $P$ is attained in that compact
set.

\bdefi With the above notations, 
\label{de:various heights}
\begin{itemize}
\item 
let $h_{M,e}$  be the infimum of all $h$ in $\RR$
such that there exists an irrational line starting from $e$ eventually
avoiding the Busemann level set $\beta_e^{-1}([h,+\infty[)$,
\item
let $h'_{M,e}$ be the infimum of the heights of the closure of
the geodesic lines contained in $C(M)$ that neither positively nor
negatively converge into a cusp, and 
\item 
let ${h''}_{M,e}$ be the infimum of the heights of the closed
geodesics.
\end{itemize}
\edefi

The {\it height spectrum} (i.e.~the subset of $\RR$ consisting of the
heights of the closed geodesic), as well as in the surface case its
restriction to simple closed geodesic, is worth more study (see for 
instance \cite{Haa,LS}).  
We give in the next section  examples were 
${h''}_{M,e}$ is attained, and examples where it is not attained (Proposition
\ref{prop:non_atteint}).


\medskip
The following result relates the Hurwitz constant to heights of closed
and non-cusp converging geodesic lines. Partial cases of the second
equality were known (see the work of Humbert and Ford (see
\cite{For2}) for SL$_2(\ZZ)$ and SL$_2(\ZZ[i])$, and \cite{HS} in the
case that $\widetilde{M}$ is isometric to $\HH^2$).

\btheo
\label{theo:equality_four_terms} 
Let $M$ be a non elementary geometrically finite pinched negatively
curved Riemannian manifold, and let $e$ be a cusp of $M$. Then
$$\frac{1}{2K_{M,e}}= \exp h_{M,e}=\exp h'_{M,e}= \exp {h''}_{M,e}.$$
\etheo

\dem 
Choose (arbitrarily) one of the parabolic fixed points in
$\partial\widetilde{M}$ corresponding to the cusp $e$, and call it
$\infty$. The other parabolic fixed points which project to $e$ are of
the form $\gamma(\infty)$ for $\gamma\in\Gamma$. For $h$ in $\RR$ and
$\gamma\in\Gamma$, let $H_{\gamma(\infty)}(h)=\gamma H_\infty(h)$ be
the horosphere centered at $\gamma(\infty)$ which is a lift of the
level set $\beta_e^{-1}(h)$. Let $\widetilde{r}$ be the geodesic
line from $\infty$ to $\gamma(\infty)$, and $r$ be the rational line,
which is the projection of $\widetilde{r}$ in the link of $e$ in
$C(M)$. By the definition \ref{def:depth} of the depth,
the (signed) distance between the intersection points of
$\widetilde{r}$ with respectively $H_{\infty}(1)$ and
$H_{\gamma(\infty)}(h)$ is $D(r) +h$.  If $\xi$ is a point in
$Lk(e,C(M))$ close enough to $r$, let $L_\xi$ be the (unique) geodesic
line starting from $\infty$ which is the closest to $\widetilde{r}$ of
the lifts of $\xi$ in $\widetilde{M}$ starting from $\infty$.  Then,
by the definition \ref{def:dist_on_the_link} of the distance
$d$ on the link, the geodesic line $L_\xi$ meets
$H_{\gamma(\infty)}(h)$ if and only if $d(r,\xi) \leq
\frac{1}{2}e^{-D(r)-h}= \frac{1}{2e^h}e^{-D(r)}$. The first
equality follows.  

\medskip
Let $H^+_{\gamma(\infty)}(h)$ be the horoball in $\widetilde{M}$ whose
boundary is $H_{\gamma(\infty)}(h)$. Let $\widetilde{C(M)}$ be the
convex hull of the limit set $\Lambda(\Gamma)$. It is immediate that
$h'_{M,e}$ is the infimum of all $h$ in $\RR$ such that
$\widetilde{C(M)}- \bigcup_{\gamma\in\Gamma} H^+_{\gamma(\infty)}(h)$
contains a geodesic line. Also, ${h''}_{M,e}$ is the infimum of all
$h$ such that $\widetilde{C(M)}- \bigcup_{\gamma\in\Gamma}
H^+_{\gamma(\infty)}(h)$ contains a periodic geodesic line (i.e.~one
whose stabilizer acts cocompactly on it).

The inequality $h'_{M,e}\leq {h''}_{M,e}$ is clear. Since $M$ is
geometrically finite and non elementary, there exists at least one
closed geodesic in $M$. Hence $h''_{M,e}$ is different $+\infty$, and so is
$h'_{M,e}$.

Let us prove that $h_{M,e}\leq h'_{M,e}$. Assume first that $h'_{M,e}$ is not $-\infty$. For every $\epsilon>0$, let
$h\in \RR$ such that $$h'_{M,e}< h \leq h'_{M,e}+\epsilon.$$ By the
definition of $h'_{M,e}$, there exists a geodesic line $c$ in
$\widetilde{C(M)}$ avoiding $\bigcup_{\gamma\in\Gamma}
H^+_{\gamma(\infty)}(h)$, so in particular its endpoints are not in
the orbit of $\infty$. Let $L_\xi$ be the geodesic line starting at
$\infty$ and ending at one of the endpoints of $c$. The geodesic line
$L_\xi$ is contained in $\widetilde{C(M)}$. By the existence of a
negative upperbound on the sectional curvature, $L_\xi$ and $c$ are
asymptotic at their common endpoint.  That is, there exists a positive
subray of $L_\xi$ which is contained in an $\epsilon$-neighborhood of
$c$. In particular this subray does not meet
$\bigcup_{\gamma\in\Gamma} H^+_{\gamma(\infty)}(h+\epsilon)$, since
the horospheres $H_{\gamma(\infty)}(h)$ and
$H_{\gamma(\infty)}(h+\epsilon)$ are at distance $\epsilon$ one from
the other, with $H^+_{\gamma(\infty)}(h+\epsilon)\subset
H^+_{\gamma(\infty)}(h)$. Hence $L_\xi$ projects into $M$ to an
irrational line meeting only finitely many projections of horoballs
$H^+_{\gamma(\infty)}(h+\epsilon)$. Therefore $h_{M,e}\leq
h+\epsilon\leq h'_{M,e}+2\epsilon$. The assertion follows. In
particular, $h_{M,e}$ is not $+\infty$. An analogous proof shows that if 
$h'_{M,e}=-\infty$, then $h_{M,e}=-\infty$.

Let us prove that ${h''}_{M,e}\leq h_{M,e}$.  Assume first that $h_{M,e}$ 
is not $-\infty$. For every $\epsilon>0$,
let $h>0$ such that $$h_{M,e}< h \leq h_{M,e}+\epsilon.$$ By the
definition of $h_{M,e}$, there exists a geodesic line $L_\xi$ starting
at $\infty$, ending in a point of $\Lambda(\Gamma)$ which is not a
parabolic fixed point, and which meets only finitely many horoballs
$H^+_{\gamma(\infty)}(h)$. Let $R$ be a positive subray of $L_\xi$
avoiding these horoballs. Since the endpoint of $R$ is a conical limit
point, the image $r$ in $M$ of the ray $R$ is recurent in a compact
subset of $M$. Moreover, it is recurent in a compact subset of
$T^1M$. Therefore $r$ comes arbitrarily close to itself in $T^1M$. By
the closing lemma (see for instance \cite{Ano}), there exists a closed geodesic
which is contained in the $\epsilon$-neighborhood of $r$. Any
pre-image of this closed geodesic is a periodic geodesic line, which
avoids every horoball $H^+_{\gamma(\infty)}(h+\epsilon)$. Hence
${h''}_{M,e}\leq h+\epsilon \leq h_{M,e}+2\epsilon$. The inequality
follows.  An analogous proof shows that if 
$h_{M,e}=-\infty$, then $h''_{M,e}=-\infty$. 
This completes the proof of the theorem.
\eop{\ref{theo:equality_four_terms}}

\bigskip
\noindent{\bf Proof of Theorems \ref{theo:existence_Hurwitz} and
\ref{theo:positivity_Hurwitz}:} Since $M$ is geometrically finite and
non elementary, there exists at least one closed geodesic in
$M$. Furthermore, every closed geodesic meets a fixed compact subset of $M$,
obtained by removing a cusp neigborhood of each end from
$C(M)$. Therefore $h''_{M,e}$ is finite.  It hence follows
that $K_{M,e}$ is positive and finite.
\eop{\ref{theo:existence_Hurwitz}, \ref{theo:positivity_Hurwitz}}

\bigskip
Till the end of this section, we assume that $M$ is a 
geometrically finite hyperbolic
$3$-orbifold, uniformized as in the end of
subsection~\ref{sect:const_Curv_Case}. For any non parabolic element
$\gamma\in\Gamma$, define the {\it height} of $\gamma$ to be
$$ht(\gamma)=e^{-\frac{D(r_\gamma)}{2}}\,
\left|\,\sinh\frac{\ell(\gamma)}{2}\,\right|$$ where $r_\gamma$ is the
rational line starting from $\infty$ corresponding to the double coset
$\Gamma_\infty\gamma\Gamma_\infty$, and $\ell(\gamma)$ is the complex
translation length of $\gamma$.

\blemm 
\label{lem:height_comput}
The height of $\gamma$ is
$$ht(\gamma)=\left|\frac{\sqrt{{\rm tr}^2 \gamma
-4}}{2c(\gamma)}\right|,$$ it is the euclidean vertical coordinate of
the highest point on the translation or rotation axis of $\gamma$.
\elemm

\dem 
Let $\gamma=\pm\left(\begin{array}{cc} a & b\\ c&
d\end{array}\right)$ be a non-parabolic element in $\Gamma$.  The
first claim follows from Lemma \ref{lem:depth_computation} and the
equality $\cosh\frac{\ell(\gamma)}{2}=\frac{{\rm tr}\,\gamma}{2}$.

The second part is well known (see \cite{HS}).  The two endpoints
$\gamma^-,\gamma^+$ of its translation or rotation axis are the
solutions of the equation $\frac{az+b}{cz+d}=z$. Hence
$\gamma^\pm=\frac{a-d\pm\sqrt{{\rm tr}^2 \gamma -4}}{2c}$ (since
$\gamma$ is not parabolic, one has $c\neq 0$ ). The highest point on the
translation axis of $\gamma$ has vertical coordinate
$\frac{1}{2}|\gamma^+-\gamma^-|$. The result follows.
\eop{\ref{lem:height_comput}}

\bcoro
\label{coro:Hurwitz_constcurv}
If $K_{M,e}$ is the Hurwitz constant of a cusp $e$ of a geometrically 
finite hyperbolic $3$-orbifold $M$, then
$$\frac{1}{2K_{M,e}}= \inf_{\{\gamma\in\Gamma\,|{\rm ~Re~}\ell(\gamma) 
>0\}}
\max_{\delta\in\Gamma} \; ht(\delta\gamma\delta^{-1}).$$
\ecoro

\dem The exponential of the height of the closed geodesic 
representing the conjugacy
class of an element $\gamma$ of $\Gamma$ is the supremum of the
Euclidean vertical coordinates of the points on the lifts to $\HH^3$
of the closed geodesic, hence is exactly $\sup_{\delta\in\Gamma} \;
ht(\delta\gamma\delta^{-1}).$ To see that this upper bound is
attained, as the translation length is a conjugacy invariant, one only
has to apply Remark \ref{rem:non_numerote} and Lemma
\ref{lem:depth_computation}, which imply that the $|c(\gamma)|$'s, for
$\gamma$ moving the point $\infty$, form a discrete subset of $\RR$
with finite multiplicities, and have a positive lower bound.
\eop{\ref{coro:Hurwitz_constcurv}}

\bigskip
\rem The right handside of the equation may seem to depend on
conjugation of $\Gamma$, but the map $ht$ has been defined by suitably
choosing some conjugate of $\Gamma$. The formula may be used to
calculate on computers the Hurwitz constants of Fuschian or Kleinian
groups, as the Hurwitz constants of the mod $p$ congruence subgroups
of SL$_2(\ZZ)$.

\medskip
\rem If $M_{\ZZ}$ is $\HH^2/{\rm PSL}_2(\ZZ)$, then $M_{\ZZ}$ has
only one end and (the first equality being due to Hurwitz)
$$\frac{1}{2K_{M_{\ZZ}}}= \exp \frac{\sqrt{5}}{2}=\exp
ht\left(\begin{array}{cc} 2 & 1 \\ 1 & 1\end{array}\right).$$ It is
proved in \cite{HS} (as well as analogous statements for the case of
Hecke groups) that $\left(\begin{array}{cc} 2 & 1 \\ 1 &
1\end{array}\right)$ realizes the maximum height in its
conjugacy class.

\medskip
If $d$ is a squarefree positive integer, if ${\cal O}_{-d}$ is the ring
of integers in the imaginary quadratic field $\QQ(\sqrt{-d})$, and
$M_{d}$ is the hyperbolic $3$-orbifold quotient of $\HH^3$ by the {\it
Bianchi group} PSL$_2({\cal O}_{-d})$, then (see for instance
\cite{Swa}) $M_d$ has one and only one cusp if and only if
$d=1,2,3,7,11,19,43,67,163$, in which case the known values and
estimates on the Hurwitz constant $K_{M_d}$
are given by the following table, up to our knowledge:

{\small
$$\begin{array}{|c|c|c|c|c|c|c|c|c|c|} 
\hline 
d & 1 & 2 & 3 & 7 & 11 & 19 &  43 & 67 & 163  \\ 
\hline  
& & & & & & & & & \\
K_{M_d} & \begin{array}{c}\frac{1}{\sqrt{3}} \\ \cite{For2}\end{array} &
 \begin{array}{c}\frac{1}{\sqrt{2}} \\ \cite{Per33}\end{array} & 
 \begin{array}{c}\frac{1}{\sqrt[4]{13}} \\ \cite{Per31}\end{array} & 
 \begin{array}{c}\frac{1}{ \sqrt[4]{8}} \\ \cite{Hof36}\end{array} &
 \begin{array}{c}\frac{2}{\sqrt{5}} \\ \cite{DP}\end{array} & 
 \begin{array}{c}1 \\ \cite{Poi}\end{array} &
\begin{array}{c}  
 \sqrt{\frac{11}{5}} \;\leq\; K_{M_{43}} \;\leq\; \sqrt{\frac{13}{3}} 
\\ 
\;\cite{Poi}\;\;\;\;\;\;\;\cite{Vul1}\end{array}  & 
? & ? \\ 
& & & &  & & & & &\\
\hline 
\end{array}$$}

The only $d$ for which we know an element of $PSL_2({\cal
O}_{-d})$ realizing the min-max in the expression of
$\frac{1}{2K_{M_d}}$ given by Corollary \ref{coro:Hurwitz_constcurv} is
$d=1$ (see \cite{For2} between the lines), for which the following
element works:
$$\left(\begin{array}{cc}2-i& 2i \\-2i & 2+i \end{array}\right)$$

Note that in particular, the formula of Corollary
\ref{coro:Hurwitz_constcurv} explains the general shape of the above
results, since in these cases, the Hurwitz constant is an inf-max of
numbers of the form $\frac{2\sqrt{N(u)}}{\sqrt[4]{N(v)}}$ with $u,v$
algebraic integers in ${\cal O}_{-d}$, hence having integral norms
$N(u), N(v)$.

\section{The Hurwitz constant of once-punctured hyperbolic tori}
\label{sect:Hurw_torus}

In this section, we study the Hurwitz constant for cusped hyperbolic
surfaces $M$ with a choosen cusp $e$, computing it in the case of
once-punctured hyperbolic tori.  Recall that ${h''}_{M,e}$ is defined
as the infimum of the heights (with respect to $e$) of the closed
geodesics in $M$

\bprop\label{prop:non_atteint} 
Let $M$ be the (unique up to isometry) complete hyperbolic thrice
punctured sphere, and $e_1, e_2, e_3$ its cusps. The infimum defining
${h''}_{M,e_1}$ is not attained on a closed geodesic, but is attained
on a (simple) geodesic line starting from $e_2$ and converging to
$e_3$.  
\eprop

\dem The manifold $M$ is obtained by doubling along its boundary an
ideal hyperbolic triangle $\tau$. Fix a lift of one of these two
triangles in a universal cover $\widetilde{M}$ of $M$, and identify it
with $\tau$ (by the covering map). Let $\ell$ be the side of $\tau$
opposite to its vertex (mapping to the cusp) $e_1$. The height of
$\ell$ is easily seen to be exactly $0$. Since $M$ has finite volume,
the pairs of endpoints of lifts of closed geodesic are dense in
$\partial\widetilde{M}\times\partial\widetilde{M}$. Therefore there exists
in $\widetilde{M}$ a sequence of translation axes of covering group
elements, whose endpoints converge to the endpoints of
$\ell$. Projecting to $M$, one obtains a sequence of closed geodesics
which converges, for the uniform Hausdorff distance on compact subsets
of $M$, to (the image in $M$ of) $\ell$. We conclude that
the heights of these closed geodesics converge to $0$. 
Consider the  open subset $U$ of points of $M$ with height
strictly less than $0$. It is the disjoint union of two open
half-cylinders whose fundamental groups are parabolic. Since 
there are no closed geodesic entirely lying in $U$, this 
proves the result.
\eop{\ref{prop:non_atteint}}

\medskip
In the same way, there exists an hyperbolic torus $M$ with two cusps
$e_+, e_-$ such that the infimum defining ${h''}_{M,e_+}$ is not
attained on a closed geodesic, but is attained on a geodesic line
starting from $e_-$ and converging to $e_-$. One may for instance take
the double along the boundary of the hyperbolic pants side lengths
$0,a,a$, i.e.~one cusp and two totally geodesic boundary of the same
length $a>0$. The difference between the infimum of the heights of
simple closed geodesic with the one on all closed curve can go to
$+\infty$, as can be seen with the previous example by letting $a$
goes to $0$.

\bigskip
Let ${\cal M}_{g,1}$ be the moduli space of one-cusped hyperbolic
metrics on a closed, connected and oriented surface $S$ with genus 
$g\geq 1$ and having one puncture.

\bprop\label{prop:minimum_Hurwitz_moduli_space} 
Let $h'':{\cal M}_{g,1}\ra ]-\infty,+\infty[$ be the map, which
associates to (the isometry class of) a one-cusped hyperbolic metric
on $S$, the infimum of the heights with respect to its cusp of its
closed geodesics. Then $h''$ is continuous and proper.  
\eprop

\dem Let $\sigma,\sigma_0$ be two one-cusped hyperbolic metric on a
punctured torus $S$. If $\sigma$ is close to $\sigma_0$ in the moduli
space, then there exists a smooth diffeomorphism $f$ of $S$ such that
$f^*\sigma,\sigma_0$ coincide outside a compact subset, and are
sufficiently close on that compact subset. In particular, their height
functions are close. Every closed geodesic $\gamma$ for one
metric is then a closed curve with geodesic curvature at most $\epsilon$
for the other metric. But such a closed curve $c$ is at uniform
distance at most $d_\epsilon$ from a genuine closed geodesic $c'$,
with $d_\epsilon$ depending only on $\epsilon$ tending to $0$ as
$\epsilon\ra 0$. Hence the height of $c$ and $c'$ are very close.
This proves that $h''$ is continuous.

Let us prove that $h''(u)$ converges to $-\infty$ when $u$ exits every
compact subset of ${\cal M}_{g,1}$. By the Mumford Lemma, hyperbolic
surfaces converge to infinity in their moduli space if and only if
they develop a short closed geodesic.  By the Margulis
Lemma, the height of a short geodesic is low, so the result
follows.  \eop{\ref{prop:minimum_Hurwitz_moduli_space}}

\medskip
This continuity and properness result holds for a much larger
class of Riemannian manifolds. For any $a,b,v$ in $]0,+\infty[$, let
$\M_{a,b,v,n}$ be the set of (isometry classes of) Riemannian
$n$-manifolds $(M,\sigma)$ with sectional curvature $K$ satisfying
$-b^2\leq K\leq -a^2$, with volume at most $v$, and with one cusp
$e$. Endow it with the Lipschitz topology on compact subsets,
i.e.~$(M,\sigma)$ is close to $(M',\sigma')$ if there exist big
compact submanifolds $K$ of $M$ and $K_0$ in $M_0$, and a smooth
diffeomorphism $f$ from $K$ to $K'$ such that $f^*\sigma,\sigma'$ are
uniformally close on $K'$.  Recall that by an easy adaptation to the finite volume with one cusp case of the compactness theorem of Cheeger-Gromov, for 
very $i>0$ the subset of points of $\M_{a,b,v,n,i}$
whose non-peripheral injectivity radius is at least $i$ is
compact. This theorem gives the analog of the Mumford lemma of
the previous situation.

\bigskip
We now give the explicit computation of the map $h'':{\cal M}_{1,1}\ra
]-\infty,+\infty[$ in terms of the Fenchel-Nielsen coordinates. In
particular we prove that $h''$ is real-analytic, that the infimum
defining $h''([\sigma])$ is attained on (one of) the shortest simple closed
geodesic for $[\sigma]\in {\cal M}_{1,1}$. We also compute the maximum of 
$h''$ on the moduli space, as well
as on which points it is attained.

We start with a few easy geometrical lemmae, whose proofs are either
omitted or sketched.

\blemm\label{lem:haut_cercle_tangent} 
In the Euclidean plane, consider two circles $c,d$ of radii $r,s$,
bounding disjoint discs, tangent to a line, with distance $t$
between the tangency points $P,Q$.  If $t\geq r+s$ and $r\geq s$, then
the radius $R$ of the half-circle $C$ orthogonal to the line, tangent to
$c,d$ and which, when starting from $[P,Q]$, meets first $d$, is bigger
than the radius $S$ of the one $D$  first meeting $c$.  
\elemm

\myfigure{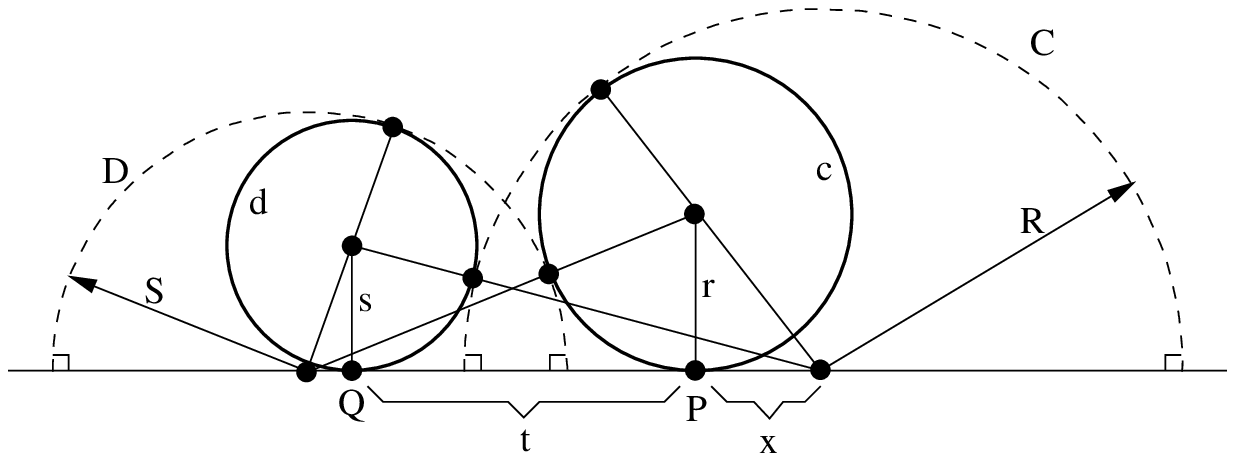}{4.7cm}{Heights of common tangent circles}

\dem By considering the right-angled triangle with vertices at the center
of $d$, the center of $C$, and $Q$, we obtain
$$(R+s)^2=s^2+ (t+x)^2,$$ with $x$ the (signed) distance between the
center of $C$ and $P$. By considering the right-angled triangle with
vertices the center of $c$, the center of $C$ and $P$, one gets
$$(R-r)^2=r^2+ x^2.$$ Eliminating $x$, this implies that
$$R^2(t^2-(r+s)^2)+Rt^2(s-r) -\frac{t^2}{4}=0.$$ Similarly, one gets
$$S^2(t^2-(r+s)^2)+St^2(r-s) -\frac{t^2}{4}=0.$$
By substracting the second equation from the first, one gets
$$(R-S)(R+S)(t^2-(s+r)^2) +(s-r)t^2(S+R)=0$$ so that if $t\geq r+s$
and $r\geq s$, then $R\geq S$.  
\eop{\ref{lem:haut_cercle_tangent}}

\blemm\label{lem:comput_hyp} 
For every $\alpha\leq\frac{\pi}{2}$, in the upper-halfspace model of
the hyperbolic plane, the distance between points of angle $\alpha$
and of angle $\frac{\pi}{2}$ on a non vertical hyperbolic geodesic is
$\log\; \cot\frac{\alpha}{2}$.  \eop{} 
\elemm

\smallskip
Let $S$ be the (smooth) once-punctured torus, and $\gamma$ an
essential (i.e.~homotopic neither to a point nor to the puncture)
simple closed curve on $S$. The real-analytic Fenchel-Nielsen 
coordinates (see for
instance \cite{FLP}) for the Teichm\"uller space ${\cal T}_{1,1}$ of
the marked hyperbolic metrics on $S$ are the length $\ell\in ]0,+\infty[$
of (the closed geodesic for the marked hyperbolic metric which is
freely homotopic to) $\gamma$ and the twist parameter $\theta\in\RR$
around it.

We define $h''(\ell,\theta)$ as the infimum of the heights of closed
geodesics on the image by the canonical map ${\cal T}_{1,1}\ra{\cal
M}_{1,1}$ of the point with Fenchel-Nielsen coordinates
$\ell\in]0,+\infty]$ and $\theta\in[-\infty,+\infty]$.

We will use the following constants. Let $\ell_{\rm min}=2\ln (1+\sqrt{2})$, 
$\ell_{\rm max}=2\ln\frac{3+\sqrt{5}}{2}$, and $$\theta_{\rm min}(\ell)=
\left\{\begin{array}{ll} \frac{4\pi}{\ell}\cosh^{-1}(\sinh\frac{\ell}{2}) &
{\rm if~} \ell\geq\ell_{\rm min}\\ 0 & {\rm ~otherwise}\end{array}\right..$$

The following is a well-known range reduction for the study of  $h''(\ell,\theta)$.

\bprop \label{prop:range_reduction}
Every point the Teichm\"uller space ${\cal T}_{1,1}$ is equivalent
under the mapping class group of $S$ to a point having Fenchel-Nielsen coordinates $(\ell,\theta)$ satisfying $\ell\in]0,\ell_{\rm max}]$, and $\theta\in[\theta_{\rm min}(\ell),\pi]$.
\eprop

\dem
Recall that the diffeomorphism group of
$S$ acts transitively on the essential simple closed curves on $S$.
Hence any marked hyperbolic structure on $S$ is equivalent under the 
mapping class group to a new one for which (the closed geodesic which is
freely homotopic to) $\gamma$ is one of the shortest closed geodesics.
Let us prove that the Fenchel-Nielsen coordinates of the new point in 
${\cal T}_{1,1}$ satisfy the above requirements.

\medskip
\noindent{\bf Step 1}: Range of $\ell$.

It is well known (see for instance \cite{Schmu}) that the unique (up to isometry) once-punctured hyperbolic torus, such that its {\it systole} 
(i.e.~the length 
of its shortest closed geodesic) is maximum, 
is the modular one $T_{\rm mod}$ (i.e.~$T_{\rm mod}=\HH^2/\Gamma$
where $\Gamma$ is the commutator subgroup of
PSL$_2(\ZZ)$). Geometrically, $T_{\rm mod}$ is obtained by gluing 
isometrically opposite faces of an hyperbolic hexagon with a dihedral 
symmetry group of order $6$, and angles alternatively $0$ and $2\pi/3$.
The closed geodesics, whose lengths are the smallest, are exactly the three 
simple geodesics obtained by taking the common
perpendicular to opposite edges of the hexagon.

\myfigure{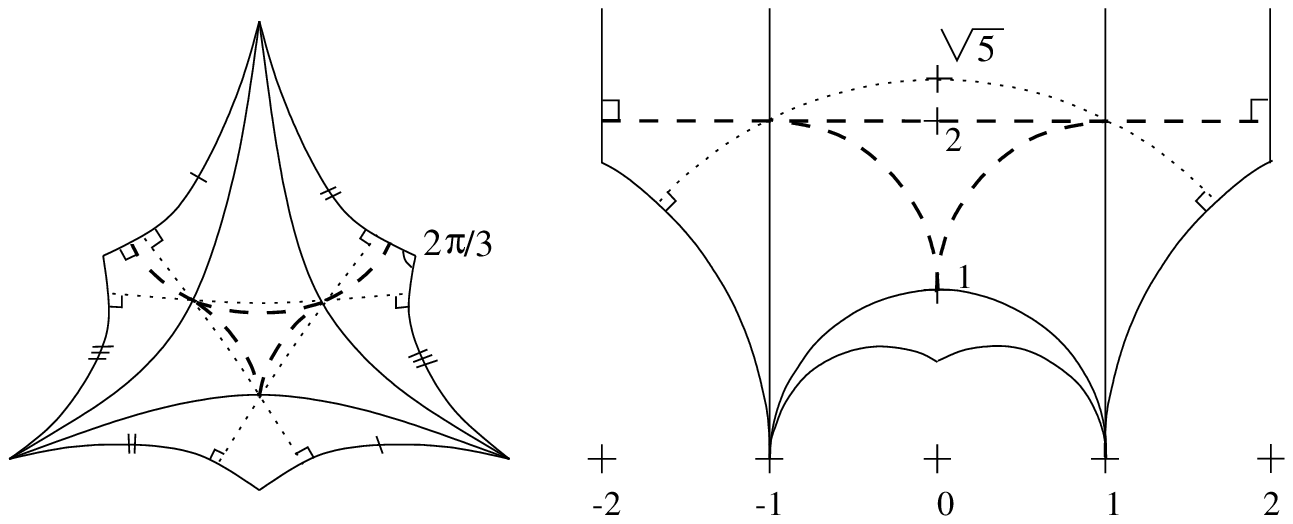}{5cm}{The modular once-punctured torus}

Using Lemma \ref{lem:comput_hyp}, it is easy to see that the Fenchel-Nielsen 
$\ell$-coordinate of $T_{\rm mod}$ (with the obvious marking) is  $$\ell=4\log\;\cot\;\arctan 2=\ell_{\rm max}.$$ 
(It follows by uniqueness and the Proposition \ref{prop:range_reduction} 
that, $\theta$ belonging to $[\theta_{\rm min}(\ell),\pi]$ and $\theta_{\rm min}(\ell_{\rm max})$ being $\pi$, the Fenchel-Nielsen 
$\theta$-coordinate of $T_{\rm mod}$ is $\pi$, but we will not need this.)

For future reference, the height of $\gamma$ is (see Figure \arabic{fig}) 
\begin{equation}\label{eq:hmax}
h''_{\rm max}=\int_{2}^{\sqrt{5}}\frac{dt}{t}=
\log\frac{\sqrt{5}}{2}.
\end{equation}

\noindent{\bf Step 2}: Range of $\theta$.

Since Dehn twists of angle multiple of $2\pi$ around $\gamma$ define 
elements of the mapping class group, and since each (complete finite 
volume) hyperbolic metric on $S$ has an elliptic involution, we
need only to consider the twist angles $\theta\in[0,\pi]$.

\medskip
Cutting open along $\gamma$, one obtains an hyperbolic pair of pants
with side lengths $(\ell,\ell,0)$. It is well-known (see for instance
\cite{FLP}) that such a pair of pants is isometric to the double of a
right-angled hyperbolic pentagon $P$ with one ideal vertex, along the
sides adjacent and opposite to the ideal vertex, the two other sides
having length $\frac{\ell}{2}$. We work in the upper halfplane model,
with the ideal point at infinity, and $P$ contained in the first quadrant,
meeting the vertical axis in $[1,+\infty[$. Let $s$ be the
highest point on the side of $P$ opposite to the vertex at infinity.

\myfigure{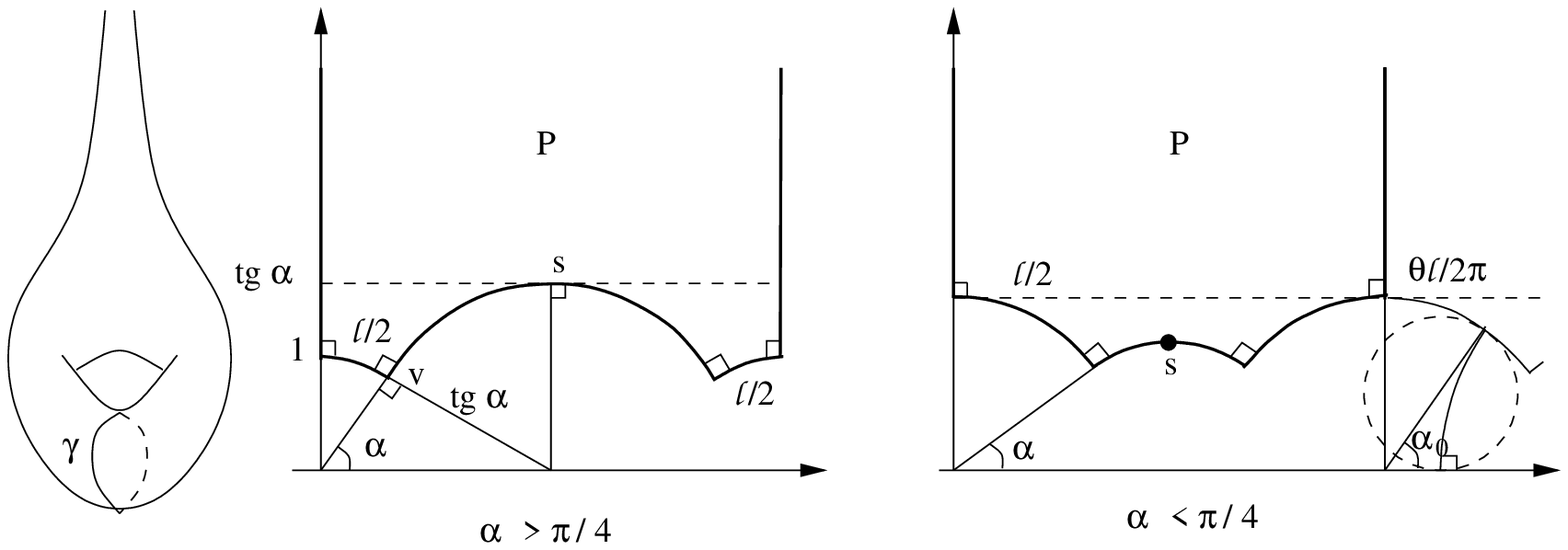}{4.6cm}{Once-punctured hyperbolic tori}

Let $u$ be the finite vertex of $P$ on the vertical axis (with vertical coordinate $1$), and $v$ be
the finite vertex adjacent to $u$. By Lemma
\ref{lem:comput_hyp}, the angle $\alpha$ of $v$ at the origin  
satisfies 
\begin{equation}\label{eq:alpha_vs_ell}
\tan \frac{\alpha}{2}=e^{-\frac{\ell}{2}}. 
\end{equation}
If $\ell$ is small enough, then the trace on $P$
of the Busemann level set $\beta_e^{-1}(0)$ is exactly the horizontal
segment through $s$. The infimum of the heights of closed geodesic is
attained exactly on $\gamma$.  A direct computation gives
$$h''(\ell,\theta)= -\log \tan \alpha = \log\;\sinh \frac{\ell}{2}.$$
This formula is valid for $\ell$ small until $s$ and $u$ are at the
same height, that is until $\ell=2\log(1+\sqrt{2})=\ell_{\rm min}$ 
(for $\alpha =\frac{\pi}{4}$). 
For future reference, we state this result as a proposition.

\bprop\label{prop:comput_Hurw_dom_low}
If $(\ell,\theta)\in\;]0,\ell_{\rm min}]\times[0,\pi]$, then 
$h''(\ell,\theta)= \log\;\sinh \frac{\ell}{2}.$ \eop{}
\eprop

Note that $h''(\ell,\theta)$ is analytic on
$(\ell,\theta)$, increasing in $\ell$ and does not depend on $\theta$ on
this range.

 The once-punctured hyperbolic torus with
$\ell=\ell_{min}$ and $\theta=0$ is the (unique up to isometry)
once-punctured hyperbolic torus with an order $4$ symmetry group,
obtained by identifying (without gliding) the opposite sides of a
regular ideal hyperbolic quadrangle. There are exactly two closed
geodesics whose heights are minimal. They are obtained by taking the
common perpendicular to the opposite sides of the quadrangle.

\medskip
As $\ell$ increases starting from $\ell_{\rm min}$, the length
$\ell'$ of the minimizing segment between the two boundary components
of the hyperbolic torus split open along $\gamma$ decreases, and
becomes shorter than the (common) length of the boundary
components. An easy computation, using Lemma \ref{lem:comput_hyp},
shows that $$\ell'=2\log\;\cot(\frac{\pi}{4}-\frac{\alpha}{2})=2\log\;
\coth \frac{\ell}{4}.$$ In order for $\gamma$ to remain (one of) the
shortest closed geodesic, we need to twist by some angle $\theta$
around $\gamma$. Let $\gamma'$ be the closed curve, obtained by
following  the path which first, in the torus split open along
$\gamma$, is the shortest common perpendicular between the two
boundary curves, and then is the shortest of the two subpaths of $\gamma$
back to its origin.  By Theorem 7.3.6 of \cite[page 183]{Bea}, recall
that the translation length $\ell(gh)$ of the product of two
hyperbolic isometries $g,h$ of the hyperbolic plane, of translation
lengths $\ell(g),\ell(h)$, with perpendicular translation axes, is
given by:
$$\cosh \frac{\ell(gh)}{2} = \cosh \frac{\ell(g)}{2}\cosh
\frac{\ell(h)}{2}$$ The closed geodesic freely homotopic to $\gamma'$
has the same length as $\gamma$ exactly when 
$$\cosh \frac{\ell}{2} = \cosh \frac{\theta\ell}{4\pi}\cosh
\frac{\ell'}{2},$$ that is when $$\theta=
\frac{4\pi}{\ell}\cosh^{-1}(\sinh\frac{\ell}{2})=\theta_{\rm min}(\ell).$$ 
Note that $\theta_{\rm min}(\ell)$ is equal to $\pi$ exactly when
$\ell$ is the length $\ell_{\rm max}$ of $\gamma$ for the modular
once-punctured hyperbolic torus. This proves the proposition about the 
range restriction.
\eop{\ref{prop:range_reduction}}

\btheo\label{theo:hauteur_maximale_moduli_space_un} 
If $\ell\in [0,\ell_{\rm max}]$ and $\theta\in [\theta_{\rm 
min}(\ell),\pi]$, then $h''(\ell,\theta)= \log \sinh \frac{\ell}{2}$.
In particular, $h''$ is real-analytic on ${\cal M}_{1,1}$. It is the 
height of one of any of the shortest closed geodesic, which is 
simple. Its maximum 
on ${\cal M}_{1,1}$ is $h_{\max}=\log\frac{\sqrt{5}}{2}$, which is 
attained uniquely on the modular one-cusped hyperbolic torus.  
\etheo

\dem 
By Proposition \ref{prop:comput_Hurw_dom_low}, we only have to 
consider the case $\ell\in [\ell_{\rm min},\ell_{\rm max}]$ and 
$\theta\in [\theta_{\rm min}(\ell),\pi]$.

\myfigure{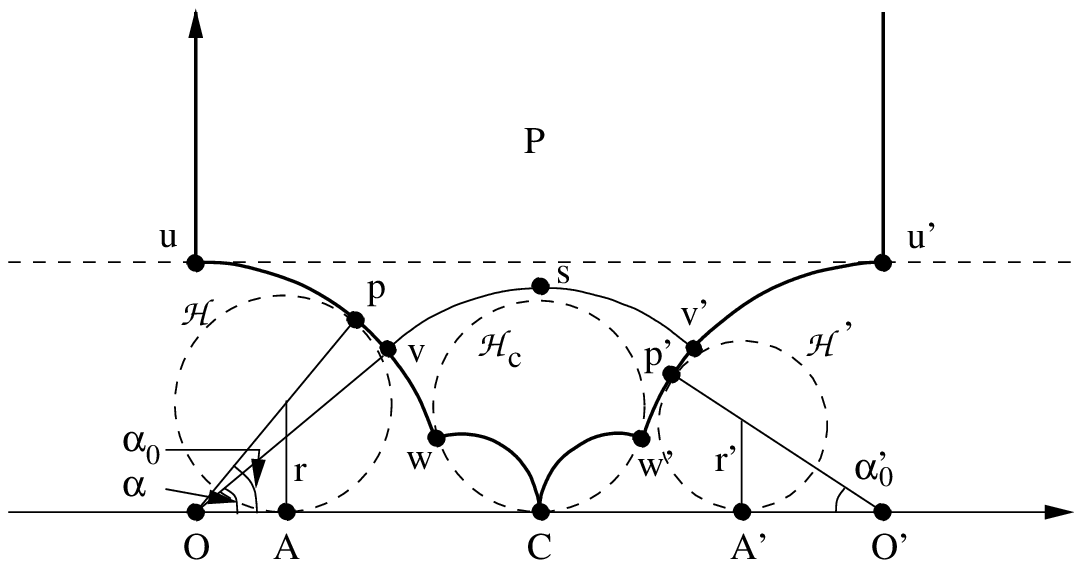}{6cm}{Computing the heights}

Take as a fundamental domain the union of the pentagon $P$ and its
image by the hyperbolic reflection along the side of $P$ opposite to 
the vertex at infinity of $P$ (see figure \arabic{fig}). Let $H_u$ be the horoball in the upper half space of points of vertical coordinates at least 
one, so that its intersection points with the boundary of the fondamental 
domain are $u,u'$. Let $\widetilde{\gamma}, \widetilde{\gamma}'$ be the 
two lifts of $\gamma$ containing  $u,u'$ respectively. 
The horosphere $\partial H_u$ is tangent to $\widetilde{\gamma}, 
\widetilde{\gamma}'$ respectively at $u,u'$. 
Define $H_c$ to be the image of $H_u$ by the reflection along the side 
of $P$ opposite to the vertex at infinity of $P$. The horosphere 
$\partial H_c$ is tangent to $\widetilde{\gamma}, \widetilde{\gamma}'$
at the vertices $w,w'$ respectively of the fondamental domain.

Recall that $v,v'$ were the vertices of $P$ opposite to the point at 
infinity, with $v$ the closer to $u$, and $\alpha$ is the angle  of $v$ at 
the origin.  Let $g$ be the element of the 
covering group that is the translation along the geodesic line through 
$v,v'$ sending $v$ to $v'$, composed with the translation along 
$\widetilde{\gamma}'$ of length $\frac{\theta\ell}{2\pi}$ (and moving 
$v'$ towards 
$u'$). Let $\H$ be the preimage of $H_u$ by $g$, and $\H'$ be the image 
of $H_c$ by $g$. Denote by $A, C,A'$ the points at infinity of the 
horoballs $\H,H_c,\H'$ respectively (see Figure \arabic{fig}).
Note that $\H$ is tangent to $\widetilde{\gamma}$ at the point $p$ at 
distance $\frac{\theta\ell}{2\pi}$ from $u$ on the side $[u,v]$ of $P$. 
Similarly $\H'$ is tangent to $\widetilde{\gamma}'$ at the point $p'$ at distance $\ell-\frac{\theta\ell}{2\pi}$ from $u'$ on 
the side $[w',u']$ of the fundamental domain.

We start with a few easy computations that will be needed in the proof.

\blemm\label{lem:diametre_Hc}
The Euclidean radius of the disc $H_c$ is 
$r_c=1/(2\sinh^2\frac{\ell}{2})$.
\elemm

\dem 
We have already seen (see the discussion after Equation (\ref{eq:alpha_vs_ell}))
that if $s$ is the highest point on the 
side of $P$ opposite to infinity, then the hyperbolic distance 
between $s$ and the horizontal horosphere $\partial H_u$ is 
$|\log \sinh\frac{\ell}{2}|$. Since $s$ is the midpoint of the 
common perpendicular segment to the horospheres $\partial H_u,\partial H_c$, 
the Euclidean diameter of the disc $H_c$ is 
$1/\sinh^2\frac{\ell}{2}$. The result follows.
\eop{\ref{lem:diametre_Hc}}

\blemm\label{lem:diametre_H}
The Euclidean radius of the disc $\H$ is 
$r= 1/(2\cosh^2\frac{\theta\ell}{4\pi})$.
\elemm

\dem 
If $\alpha_0$ is the angle of the center of the 
disc $\H$ at the origin, then by Lemma \ref{lem:comput_hyp}, one has
$\tan\frac{\alpha_0}{2} = e^{-\frac{\theta\ell}{2\pi}}$.
Since $\partial \H$ is tangent to $\widetilde{\gamma}$ at $p$,
one has $\sin\alpha_0=\frac{r}{1-r}$. The result follows. \eop{\ref{lem:diametre_H}}

\blemm\label{lem:distance_AC}
The Euclidean distance between $A$ and $C$ is $\coth \frac{\ell}{2}-\tanh\frac{\theta\ell}{4\pi}$.
\elemm

\dem 
Let $O$ be the origin in the plane. One has $d(A,C) = d(C,0)- d(0,A)$ 
which gives 
$d(A,C)=\frac{1}{\cos \alpha}-\sqrt{(1-r)^2 -r^2}$.
Using Equation (\ref{eq:alpha_vs_ell}) and Lemma \ref{lem:diametre_H},
the result follows. \eop{\ref{lem:distance_AC}}

\bigskip
\noindent{\bf Step 1}: A computation of the Busemann level $\beta_e^{-1}(0)$.

It is easy to see that the only horospheres in the orbit of $H_u$, 
that meet the fundamental domain, are $H_u,\H,\H_c,\H'$ (unless $\theta=0$,
where we have two more horospheres, the translations of $\H,\H'$ along $\widetilde{\gamma},\widetilde{\gamma}'$ of a distance $\ell$, that meet the fundamental domain in $w,u'$ respectively.)
Note that since $\theta\geq\theta_{\rm min}(\ell)$, the Euclidean
radius of $\H$, which is $r=1/(2\cosh^2\frac{\theta\ell}{4\pi})$ by 
Lemma \ref{lem:diametre_H}, is less than the Euclidean
radius of $H_c$, which is $1/(2\sinh ^2\frac{\ell}{2})$ by Lemma 
\ref{lem:diametre_Hc}.

Hence the Busemann level set $\beta_e^{-1}(0)$ has a lift which is the horizontal horosphere through $s$. In particular the height
of the closed geodesic $\gamma$ is $\log\sinh\frac{\ell}{2}$ since
no lift of $\gamma$ enters the interior of the horosphere $H_u$.

\medskip
\noindent{\bf Step 2}: Computation of the minimal height of a closed 
geodesic.

One only has to prove that there is no geodesic line $L$ 
meeting the fundamental domain and avoiding the horoballs 
$H_u,\H,H_c,\H'$. By absurd, assume that such a line $L$ exists.
The boundary of the fundamental domain, from which one removes the points 
lying in one of these horoballs, is the disjoint union of $4$ geodesic 
segments $I_1=]u,p[, I_2=]p,w[, I_3=]w',p'[, I_4=]p',u'[$.
The geodesic line enters the fundamental domain through one of $I_1,I_2$ 
and exits it through one of $I_3,I_4$.

Let us first prove that $L$ cannot enter through $I_1$ and exit 
through $I_4$. The diameter of the Euclidean halfcircle $L$ 
perpendicular to the real axis would be at least the distance 
between the tangent points of $\H$ and $\H'$ to the real line, 
plus the Euclidean radius of the discs $\H$ and $\H'$.
By Lemma \ref{lem:diametre_H} and \ref{lem:distance_AC}, one has
$$r= \frac{1}{2\cosh^2\frac{\theta\ell}{4\pi}} \;\;\;\;\;{\rm and}\;\;\;\;\;
d(A,C)=\coth\frac{\ell}{2} -\tanh\frac{\theta\ell}{4\pi}.$$
Similarly, with $r'$ the Euclidean radius of $\H'$, one gets $$r'= \frac{1}{2\cosh^2(\frac{\ell}{2}-\frac{\theta\ell}{4\pi})} \;\;\;\;\;{\rm and}\;\;\;\;\;d(A',C)= \coth\frac{\ell}{2} -\tanh(\frac{\ell}{2}-\frac{\theta\ell}{4\pi}).$$
So that, with $f(\ell,\theta)=r+d(A,A')+r'$, we have $$f(\ell,\theta)=2 \coth\frac{\ell}{2} +(\frac{1}{2\cosh^2\frac{\theta\ell}{4\pi}}+
\frac{1}{2\cosh^2(\frac{\ell}{2}-\frac{\theta\ell}{4\pi})}) -(\tanh(\frac{\theta\ell}{4\pi}+ 
\tanh\frac{\frac{\ell}{2}-\theta\ell}{4\pi})).$$

It is easy to see that  $f(\ell,\theta)$ is decreasing in $\ell$ since $1-\frac{\theta}{2\pi}>0$. Since $\cosh x= \cosh (-x)$, and since by 
taking the derivative, the map
$t\mapsto\tanh t+\tanh(x-t)$ is increasing in $t$ for $t\leq 2x$, 
the function $f(\ell,\theta)$ is also decreasing in $\theta$. Hence 
$$f(\ell,\theta)\geq f(\ell_{\rm max},\pi)=
2\coth\frac{\ell_{\rm max}}{2} + 1/\cosh^2\frac{\ell_{\rm max}}{4} - 
2 \tanh \frac{\ell_{\rm max}}{4}$$ 
which is about $2.58885438$, hence strictly more than $2$.
In particular, the Euclidean radius of $L$ would be strictly more than 
one, which contradicts the fact that $L$ does no enter $H_u$.

\medskip
Let us prove that $L$ cannot enter through $I_2$. One only has to show 
that an Euclidean halfcircle $L$ centered on the real axis, bounding an 
halfdisc $D$ that contains $H_c$ and does not contain $\H$ has radius at 
least one. The Euclidean radius of $L$ is at least as big as 
the one of the halfcircle
$L'$ which is tangent to both $\partial H_c$ and $\partial \H$, starts 
from the segment between the tangency points $A,C$ and first meet 
$\partial \H$ and then $\partial H_c$. But recall that $\widetilde{\gamma}$
is tangent to both $\partial \H,\partial H_c$ and has radius $1$,
and that the Euclidean radius of $H_c$ is at least the 
Euclidean radius of $\H$. By Lemma \ref{lem:haut_cercle_tangent},
to prove that $L'$ hence $L$ has radius at least $1$, one only has to 
prove that $t=d(A,C)\geq r+r_c$ is positive.

By Lemmae \ref{lem:diametre_Hc}, \ref{lem:diametre_H} and 
\ref{lem:distance_AC}, we have 
$$t=\coth \frac{\ell}{2}-\tanh\frac{\theta\ell}{4\pi}-  
\frac{1}{2\cosh^2\frac{\theta\ell}{4\pi}}-\frac{1}{2\sinh^2\frac{\ell}{2}}
= \frac{\sinh(\ell) -1}{2\sinh^2\frac{\ell}{2}} -
\frac{\sinh\frac{\theta\ell}{2\pi} +1}{2\cosh^2\frac{\theta\ell}{4\pi}}.$$
This is an increasing function in $\theta$, hence its values are
greater than or equal to its value at $\theta=\theta_{\rm min}(\ell)$.
So $$t\geq 
\frac{\sinh(\ell) -1}{2\sinh^2\frac{\ell}{2}} -
\frac{2\sinh\frac{\ell}{2}\;\sqrt{\sinh^2\frac{\ell}{2}-1} +1}{2\sinh^2\frac{\ell}{2}}=
\frac{\sinh\frac{\ell}{2}(\cosh\frac{\ell}{2}-
\sqrt{\cosh^2\frac{\ell}{2}-2}\,)-1}{\sinh^2\frac{\ell}{2}}.$$
The numerator is, by an easy derivative computation, a decreasing function 
of $\ell$ on $[\ell_{\rm min},\ell_{\rm max}]$, whose value at 
$\ell=\ell_{\rm max}$ is about $0.118$, hence is positive.
This proves that $t$ is indeed positive.

Similarly, one proves that $L$ cannot exit through $I_3$, which proves the
claim. Theorem \ref{theo:hauteur_maximale_moduli_space_un} now follows. The  maximum of $h''$ been reached uniquely when $\ell=\ell_{\rm max},\theta=\pi$, that is for the modular once-punctured torus, its value $h''_{\rm max}$ has 
been computed in Equation (\ref{eq:hmax}). This ends the proof.
\eop{\ref{theo:hauteur_maximale_moduli_space_un}}

\medskip
Theorem \ref{theo:Hurw_intro} of the introduction follows from Theorem \ref{theo:hauteur_maximale_moduli_space_un}.

\section{The cut locus of a cusp}
\label{sect:cut_locus}

We keep the notations of the beginning of section
\ref{sect:rational_rays}.  The following definition is due to
\cite[section 4]{EP} in the case that $\widetilde{M}=\HH^{n}_{\RR}$
and $M$ has one cusp.

\bdefi ({\bf Cut locus of a cusp}).  
The {\it cut locus} of the cusp $e$ in $M$ is the subset
$\Sigma=\Sigma(e)$ of points $x$ in $M$ from which start at least two
(globally) minimizing geodesic rays converging to $e$.  
\edefi

(In the case $M$ is an hyperbolic $3$-orbifold and $x$ is a singular
point, we have to count the geodesic rays with multiplicities.)

The definition implies that for any $\eta>0$, such geodesic rays meet
perpendicularly the level set $\beta_e^{-1}(\eta)$ in one and only one
point, and the lengths of their subsegment between $x$ and
$\beta_e^{-1}(\eta)$ are equal.

There is a {\it canonical retraction} from $M-\Sigma(e)$ into
$\beta_e^{-1}(1)$, which associates to a point $x$ in
$M-\Sigma(e)$ the unique intersection point with $\beta_e^{-1}(1)$ of the
unique minimizing geodesic ray starting from $x$ and converging to
$e$, if $x$ does not lie in $\beta_e^{-1}([1,+\infty[)$, or the obvious projection to $\beta_e^{-1}(1)$ otherwise. 
 (Note that this retraction is a strong deformation retract).

\medskip
Let $N:M\ra \NN-\{0\}$ be the map which assigns to each $x\in M$ the
number $N(x)$ of minimizing geodesic rays starting from $x$ and
converging to $e$. (This number is finite since the curvature is
negative.) The cut locus of $e$ is by definition
$N^{-1}([2,+\infty[)$. It is immediate that $N$ is upper
semicontinuous.

In particular, $\Sigma$ is closed. The set $M$ has a natural partition
by the connected components of $N^{-1}(\{k\})$ for $k$ in
$\NN-\{0\}$. By the same proof as in \cite[Theorem A]{Sug}, it is easy
to show that $N^{-1}(2)$ is a codimension one submanifold of $M$,
which is open and dense in $\Sigma$. We will denote $N^{-1}(2)$ by
$\Sigma_0$.

\medskip
Let $H$ and $H'$ be horospheres in $\widetilde{M}$ whose horoballs are
disjoint, the {\it equidistant subspace} of $H$ and $H'$ is by
definition, the set of points in $\widetilde{M}$ which are at the same
distance from $H$ and $H'$.

Let $x\in \Sigma_0$ and $r_1,r_2$ be the two minimizing rays starting
from $x$ and converging to $e$. Let $\widetilde{x}$ be a lift of $x$
in $\widetilde{M}$, and $\widetilde{r_1},\widetilde{r_2}$ be the lifts
of $r_1, r_2$ starting from $x$. Let $H_1,H_2$ be the (disjoint)
horospheres centered at the points
$\widetilde{r_1}(\infty),\widetilde{r_2}(\infty)$ respectively, and
covering $\beta_e^{-1}(1)$. The component of $\Sigma_0$ which contains
$x$ is the image by the covering map $\widetilde{M}\ra M$ of an open
subset of the equidistant subspace of $H_1$ and $H_2$. 

In particular, if the curvature is constant, the equidistant subspace of 
two horospheres $H,H'$ bounding disjoint horoballs is a hyperbolic 
hyperplane, hence is totally geodesic.
So that each component of $\Sigma_0$ is (locally) totally geodesic.
Furthermore, the equidistant subspace
 is the unique hyperbolic hyperplane orthogonal to the geodesic line
$L$ between the points at infinity of $H,H'$ and passing through the
point of $L$ which is equidistant from $H$ and $H'$. 
Using the transitivity
properties of Isom$_+(\HH^{n}_{\RR})$ this can be easily seen.


\medskip
A {\it stratification} of a smooth manifold $M$ is a partition of $M$
into connected smooth submanifolds called {\it strata}, such that the
closure of each stratum locally meets only finitely many strata. (See
\cite{Tro1} for a general survey about topological
stratifications. R.~Thom's definition (see \cite{Tho} or \cite[page
234]{Tro1}) required only the existence of finitely many strata in the 
closure of any
stratum, but a local such assumption is sufficient when dealing with
local properties.)

Let $X,Y$ with $Y$ contained in the closure of $X$ be any strata.  The
stratification is called {\it (a)-regular} (in the Whitney's sense,
see \cite{Whi}, with applications to analytic varieties, or
\cite[page 235]{Tro1}) if for every $x_i\in X$ converging to $y\in Y$
such that the tangent subspace $T_{x_i}X$ converges to some tangent
subspace $\tau$, it follows that $T_yY$ is contained in $\tau$.

This condition is precisely the one needed for the transversality of a
submanifold to the stratification to be stable in the smooth topology
(see \cite{Fel} or \cite[page 237]{Tro1}).  We now define the
technical assumptions refered to in the introduction that we will need
on the metric.

\bdefi {\rm (}{\bf Cute cut locus}{\rm )}.
\label{defi:cute}
The cut locus $\Sigma$ of the cusp $e$ will be called {\it cute} if
\begin{itemize}
\item 
the partition of $M$ by components of $N^{-1}(x), x\in M$
is a locally finite (a)-regular stratification of $M$.
\item
each component $\sigma$ of $\Sigma_0$ is simply connected, and has a 
unique locally highest point $\widehat{\sigma}$, that belongs to $\sigma$.
\item
  $\Sigma_0$ has only finitely many components.
\end{itemize}
\edefi

The union of codimension $0$ strata is exactly $N^{-1}(2)=\Sigma_0$.
Since any component of $\Sigma_0$ is simply connected, it has a
transverse orientation, uniquely determined by the orientation of any
transversal subspace to the tangent subspace at any point.

It follows from the fact that the highest point $\widehat{\sigma}$ is
unique and in the interior of $\sigma$ that $\widehat{\sigma}$ is the
unique point of $\sigma$ from which start perpendicularly to $\sigma$
at least two (and indeed exactly two) minimizing geodesic rays
converging to $e$.

We will call $\widehat{\sigma}$ the {\it summit} of $\sigma$. A {\it
summit} of $\Sigma$ is the summit of some component of $\Sigma_0$.  By
definition, the cut locus $\Sigma$ of $e$ does not meet the
neighborhood $\beta_e^{-1}(]0,+\infty[)$ of $e$.

Since $M$ is geometrically finite, if there is only one cusp, if it
satisfies the first condition of Definition \ref{defi:cute}, 
then the number of components of $\Sigma_0$
is finite (that is the third condition of Definition \ref{defi:cute} 
is automatically satisfied).

\bdefi {\rm (}{\bf Integral geodesic}{\rm )}.
\label{def:integral}
A geodesic line starting from $e$ whose first hitting point with 
$\Sigma$ is a summit of $\Sigma$ will be called an {\it integral line}.  
\edefi

Note that an integral line is rational, since the union of two
minimizing geodesic rays starting from a summit $\widehat{\sigma}$ of
$\Sigma$ gives by the remark following Definition \ref{defi:cute} (the
image of) a geodesic line converging both positively and negatively to
$e$. Since $\Sigma_0$ has only finitely many components, and since
there are two integral lines per summit of $\Sigma$, the set
${\cal R}={\cal R}(e)$ of rational rays is a finite subset of
$Lk(e,C(M))$.

\medskip
Let us describe interesting examples when the cut locus of the cusp is cute.
This should also be the case for other rank one symmetric spaces of non 
compact type.

\bprop\label{prop:cute_exist}
(i)~ If the curvature is constant, then $\Sigma$ is cute.

(ii)~ If the curvature is constant, if we have only one end,
 and if the codimension of each
component of $N^{-1}(\{k\})$ is $k-1$, then for any small enough
perturbation with compact support of the metric (in the C$^\infty$ 
topology), the cut locus of the cusp is cute.  
\eprop

\dem (i)~Since being a locally finite (a)-regular stratification is a
local property, one can work in the universal cover $\widetilde{M}$ of
$M$. Take the projective Klein model for $\widetilde{M}$. Consider the
union $A$ of all equidistant subspaces of pairs of horospheres
covering $\beta_e^{-1}(1)$. This is a locally
finite union of analytic submanifolds (linear ones), hence the
dimension stratification of Whitney is a locally finite (a)-regular
stratification (see \cite{Whi}). Since the lift $\widetilde{\Sigma}$
of $\Sigma$ is a closed saturated subset of $A$, the same thing holds for
$\widetilde{\Sigma}$.

Let $x$ be a point in a component $\sigma$ of $\Sigma_0$, and
$r_1,r_2$ the two minimizing rays starting from $x$ and converging to
$e$. Let $\widetilde{x}$ be a lift of $x$ in $\widetilde{M}$, and
$\widetilde{r_1},\widetilde{r_2}$ be the lifts of $r_1, r_2$ starting
from $x$. Let $H_1,H_2$ be the (disjoint) horospheres centered at the
points at infinity $a_1,a_2$ of $\widetilde{r_1},\widetilde{r_2}$ and
covering $\beta_e^{-1}(1)$.  Let $z$ be the intersection point of the
equidistant subspace of $H_1,H_2$ with the geodesic between $a_1,a_2$.

Let us prove the claim that $z$ belongs to the lift
$\widetilde{\sigma}$ of $\sigma$ containing $\widetilde{x}$, and its
image in $M$ is the unique summit of $\sigma$. 

Let $\ell$ be the Busemann distance between $\widetilde{x}$ and $a_1$,
and $m$ be the Busemann distance between $z$ and $a_1$. Let $a_3$ in
$\partial\widetilde{M}$ be any point mapping onto $e$ other than
$a_1,a_2$, and $\ell'$ (resp.~$m'$) be the Busemann distance between
$\widetilde{x}$ and $a_3$ (resp.~$z$ and $a_3$). Note that by
assumption $\ell'>\ell$. Since the hyperbolic triangle with vertices
$\widetilde{x},z,a_1$ is rectangle at $z$, it follows from hyperbolic
triangle formulae that $m'>m$. Hence $z$ belongs to the lift of
$\Sigma_0$ and is hence a summit.  By convexity, the geodesic segment
between $z$ and $\widetilde{x}$ is contained in the lift of
$\Sigma_0$. Since the Buseman function is strictly increasing along the 
geodesic segment from $\widetilde{x}$ to $z$, it follows that $x$  
cannot be a summit. Hence the claim is proved.

Let us prove now that $\sigma$ is simply connected. By absurd,
assume that there is a closed path $c$ based at $x$ which is not homotopic 
to $0$ in $\sigma$. The closure of $\sigma$ is locally convex, by convexity 
of the Busemann functions. Hence the lift of $c$ starting at 
$\widetilde{x}$ ends in $\gamma\widetilde{x}$ for some non trivial element
$\gamma$ in the covering group. By analyticity, $\Gamma$ preserves the
equidistant subspace $E$ of $H_1,H_2$. 
 By convexity, the geodesic between $\widetilde{x}$ and $\gamma\widetilde{x}$
is contained in $\widetilde{\sigma}$, and by the previous claim, so is 
the geodesic segment between $z$ and $\gamma z$. But this is a 
contradiction, since the midpoint $m$ of $[z,\gamma z]$ is contained in 
the equidistant plane of the horospheres $H_1$ and $\gamma H_1$.

Now the third property of Definition \ref{defi:cute} follows from the fact 
that the parabolic points are bounded, and that the equidistant subspaces 
are (in the Klein model of the hyperbolic space) affine subspaces.

(ii)~Under the assumptions of (ii), the third property is automatically satisfied. The first property of a cute cut
locus follows by transversality arguments. Since the cut locus in 
constant negative curvature is locally finite, and since the deformation 
is supported on a compact subset, if the deformation is small enough,
the components of $\Sigma_0$ will remain simply connected, and by strict 
convexity of the horospheres, there will remain one and only one new summit
close to each old summit, and no new one will be created far from the old ones.
\eop{\ref{prop:cute_exist}}

\medskip
The assumption in (ii) that the codimension of $N^{-1}(\{k\})$ is $k-1$
cannot be omitted, since if some cut locus in dimension 3 has a vertex of 
degree 4, then by small perturbation of the metric, the cut locus can be 
made not locally finite, as B.~Bowditch showed us.

\bigskip
Till the end of this section, we assume that $M$ is a non elementary
geometrically finite hyperbolic $3$-orbifold, uniformized as in
subsection \ref{sect:const_Curv_Case}. All the statements below extend
easily to $\HH^{n}_{\RR}$.

It follows from Proposition \ref{prop:cute_exist}, and the fact that
there cannot be any summit close to another cusp, that $\Sigma$ is a
finite piecewise hyperbolic cell $2$-complex, whose $2$-cells are
compact or finite volume hyperbolic polygons, with one vertex at
infinity for each end besides $e$.

The {\it summit} of an edge $\tau$ of $\Sigma$ is the unique point
$\widehat{\tau}$ (contained in the interior of the edge) from which
start perpendicularly to the edge at least three minimizing geodesic
segments ending in $e$.

\medskip
The following is another description of $\Sigma$. It is due to Ford
\cite{For} for some arithmetical cases, extended by Swan \cite{Swa}
for all Bianchi groups, and to \cite[section 4]{EP} when
$\widetilde{M}=\HH^{n}_{\RR}$ and $M$ has one cusp.

Let $\gamma\in\Gamma-\Gamma_\infty$, the {\it isometric sphere}
$S_\gamma$ of $\gamma$ is the hyperbolic plane in $\HH^3$, defined as
the intersection of the upper halfspace with the Euclidean sphere
centered at $\gamma^{-1}(\infty)$ and of radius
$\frac{1}{|c(\gamma)|}$.  It follows that $\gamma$ maps the horizontal
horosphere $H_\infty(h)$ defined by the equation $t=h$ centered at $\infty$ to
the horosphere $H_{\gamma(\infty)}(h)$ centered at $\gamma(\infty)$,
which is the Euclidean sphere in the upper halfspace, tangent to $\CC$
at $\gamma(\infty)$ and with diameter $\frac{1}{h|c(\gamma)|^2}$.  We
denote by $H^+_{\gamma(\infty)}(h)$ the horoball bounded by
$H_{\gamma(\infty)}(h)$.

The key point is the following fact.

\blemm 
\label{lem:key_isometric_sphere} 
For any $\gamma\in\Gamma-\Gamma_\infty$ and $h>0$ big enough, the 
isometric sphere $S_\gamma$ is the equidistant subspace of the 
horospheres $H_\infty(h)$ and $H_{\gamma^{-1}(\infty)}(h)$.  
\elemm

\dem Let $x,y$ be the points on a vertical geodesic line $L$ at
vertical coordinates $t=h$ and $t=\frac{1}{h|c(\gamma^{-1})|^2}$
respectively.  Then the midpoint of $[x,y]$ is the point on $L$ at
height $t=\frac{1}{|c(\gamma^{-1})|}$ for any $h>0$.  The hyperbolic
plane $S_\gamma$ is orthogonal to the (vertical) geodesic line between
$\infty$ and $\gamma^{-1}(\infty)$. Hence, by uniqueness, $S_\gamma$
is the equidistant subspace of $H_\infty(h)$ and
$H_{\gamma^{-1}(\infty)}(h)$.  
\eop{\ref{lem:key_isometric_sphere}}

\bigskip 
Let $S_\gamma^-$ be the half ball bounded by $S_\gamma$. It is easy to
see that $\gamma$ maps $ \HH^3- S_\gamma^-$ to $S_{\gamma^{-1}}^-$.
Define $$B_\infty = \HH^3- \bigcup_{\gamma\in\Gamma-\Gamma_\infty}
S_\gamma^-,$$ that we will call the {\it basin of center $\infty$}.
It follows that 
\begin{enumerate}
\item 
$B_\infty$ is invariant by $\Gamma_\infty$, 
\item 
no element of $\Gamma-\Gamma_\infty$ maps an interior point of
$B_\infty$ to another interior point of $B_\infty$, and 
\item  
the images of $B_\infty$ by $\Gamma$ cover $\HH^3$.
\end{enumerate}
Therefore $B_\infty/\Gamma_\infty$ is an orbifold whose boundary is a
piecewise hyperbolic
$2$-orbifold, and $M$ is obtained by pairing faces of
$B_\infty/\Gamma_\infty$. Since $M$ is geometrically finite, the
boundary of $B_\infty$ has finitely many $2$-cells up to the action of
$\Gamma_\infty$. For $\gamma\in\Gamma$, define the {\it basin of
center $\gamma(\infty)$} to be $B_{\gamma(\infty)}=\gamma B_\infty$.

\myfigure{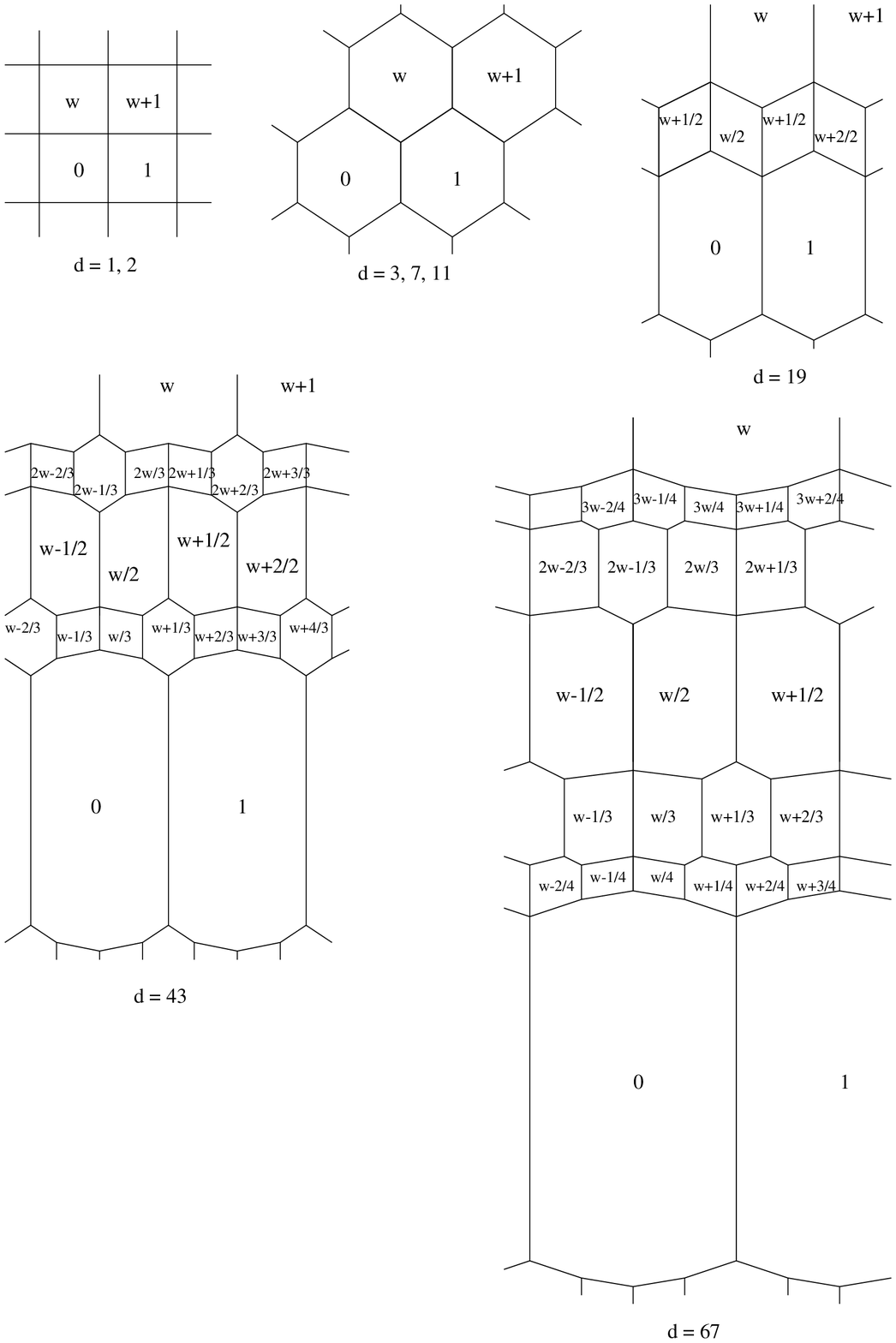}{20.5cm}{Cell decomposition of $\partial B_\infty$
 for $\HH^3/\,{\rm PSL}_2({\cal O}_{-d})$}

Above is a combinatorial picture of the vertical projection to $\CC$
of the cell decomposition of the boundary of the basin at infinity for 
$\HH^3/\,{\rm PSL}_2({\cal O}_{-d})$. Up to $d=43$, these drawings are 
due to \cite{Hat}. The complex number $w$ is $i\sqrt{d}$ if 
$d\not\equiv -1 [{\rm mod~} 4]$ and $\frac{1}{2}(1+i\sqrt{d})$ otherwise.

\medskip
Using Lemma \ref{lem:key_isometric_sphere}, the points in $\partial
B_\infty$ are exactly the points from which starts a geodesic ray
converging to $\infty$ and a geodesic ray converging to a point
different form $\infty$ in the orbit of $\infty$ by $\Gamma$. Hence
the cut locus $\Sigma$ is the image of $\partial B_\infty$ by the
canonical projection $\HH^3\ra \HH^3/\Gamma =M$. Furthermore,
the summits of the $2$-cells of $\Sigma$ are the images in $M$ of the
highest points of the isometric spheres $S_\gamma$ containing a
$2$-cell of $\partial B_\infty$, and the summits of the edges of
$\Sigma$ are the image in $M$ of highest points of the halfcircle
geodesic lines containing an edge of $\partial B_\infty$. Let
$\widetilde{{\cal R}}$ be the set of endpoints of the lifts to $\HH^3$
starting from $\infty$ of the integral lines. If $P$ is the set of
elements in $\Gamma-\Gamma_\infty$ whose isometric sphere carries a
$2$-cell of the boundary of the basin at infinity, then
$\widetilde{{\cal R}}$ is exactly the discrete subset $P(\infty)$.

In the case of the Bianchi orbifolds PSL$_2({\cal
O}_{-d})\backslash \HH^3$, the uniformization chosen in section
\ref{sect:rational_rays} is precisely the natural one 
$\HH^3\ra \HH^3/$PSL$_2({\cal O}_{-d})$. Though $\widetilde{{\cal R}}$
consists of the integer points in the case $d=1,2,3,7,11$ for
instance, it contains in general rational points with finitely many
possible denominators ($q= \pm 1, \pm 2, \pm 3$ for $d=43$, for
instance, see \cite{Poi} for other lists), not all rational lines with
these denominators being in $\widetilde{{\cal R}}$ (see Figure
\arabic{fig}).

\section{Good approximating sequence of an irrational line}
\label{sect_good_approximation}

We keep the same notations as in the beginning of section
\ref{sect:rational_rays}. In section \ref{sect:Hurwitz}, we proved
that any irrational line starting from a given cusp can be well
approximated by rational lines. The aim of this section is to give an
explicit sequence of such rational lines.

Assume that the cut locus of the cusp $e$ is cute. We say that a
geodesic line starting from $e$ is {\it totally irrational} if it recurrent 
in a compact subset of $M$ (in particular and equivalently in the finite volume case,  converges into no cusp), and is transverse to the stratification
of the cusp, i.e.~it does not meet the singular locus $\Sigma-\Sigma_0$
and is transverse to $\Sigma_0$.  This last assumption is implied by
the previous one in the case of constant curvature, since the components
of $\Sigma_0$ are (locally) totally geodesic. The set of totally
irrational lines is a dense $G_\delta$ set in the link of $e$.


Consider a totally irrational line $\xi$. Note that $\Sigma$ is cute and 
$\xi$ is transverse to $\Sigma$, hence the intersection of $\Sigma$ and 
$\xi$ has no accumulation point. Since $\xi$ is recurrent in a compact 
subset, for every $\epsilon>0$, the $\epsilon$-neighborhood of any 
positive subray of $\xi$ contains a closed geodesic (by the closing 
lemma), hence since no closed geodesic is contained in $\M-\Sigma$ 
(which retracts onto $\beta_e^{-1}(1)$), there are infinitely many 
intersection points of $\xi$ and $\Sigma$. Let $(x_n)_{n\in \NN}$ be
the sequence of intersection points of $\xi$ with
$(\sigma_n)_{n\in\NN}$, the sequence of connected components of
$\Sigma_0$ consecutively passed through by $\xi$. For each $n\in \NN$,
let $c_n$ be the path consisting of the subsegment of $\xi$ between
$e$ and $x_n$, followed by the minimizing geodesic ray from $x_n$ to
the cusp $e$, starting on the other side of $\Sigma_n$ than the one on
which $\xi$ arrives at $x_n$.  Let $r_n$ be the unique geodesic line
starting from $e$, properly homotopic to $c_n$. Since $r_n$ converges
to $e$, it is a rational line.

\bdefi
The sequence $(r_n)_{n\in \NN}$ will be called the good approximating
sequence of $\xi$, and $r_n$ the $n$-th good approximant.
\edefi

This definition can be extended to the other irrational lines, since by 
the (a)-regularity of the stratification of the cut locus, the 
transversality to the stratification is stable, at the expense of 
loosing the uniqueness of the sequence, and allowing, for rational 
lines and lines ending in other cusps, each sequence to be finite.

\medskip
\rem When $M=\HH^2/{\rm PSL}_2(\ZZ)$, this sequence coincides with the
usual best approximation sequence in the following sense: let
$\widetilde{\xi}$ be the endpoint on the real axis of the lift of
$\xi$ starting from $\infty$ ending in $[0,1[$. Note that $\xi$ is
totally irrational if and only if $\widetilde{\xi}$ is an irrational
real number. Then the endpoint on the real axis of the lift of $r_n$
starting from $\infty$ and ending in $[0,1[$ is $p_n/q_n$, where
$p_n/q_n$ is the $n$-th convergent of the irrational real number
$\widetilde{\xi}$. In the case of $M=\HH^3/{\rm PSL}_2({\cal
O}_{-d})$, we don't know if our good approximation sequence coincides
with Poitou's best approximation sequence in \cite{Poi}.

\medskip
The good approximating sequence nicely approximates an irrational line.
We prove this for totally irrational lines, but the result could be
extended for general irrational lines.

\btheo\label{theo:good_indeed}
Let $M$ be a non elementary geometrically finite pinched negatively
curved Riemannian manifold, with only one cusp $e$, having a cute cut
locus. There exists a constant $c>0$ such that, for every totally
irrational line $\xi$ starting from $e$ with good approximating
sequence $(r_n)_{n\in\NN}$, for every $n$ in $\NN$,
$$d(r_n,\xi)\leq  c\;e^{-D(r_n)}.$$
\etheo

\dem 
Let $$c=\frac{1}{2}\;\max_{x\in\Sigma} e^{-\beta_e(x)}.$$ Since
$\Sigma$ is in the complement of $\beta_e^{-1}(]0,+\infty[)$,
$\sup_{x\in\Sigma} -\beta_e(x)$ is strictly bigger than $0$. The
supremum is a maximum, since $\Sigma$ is non empty and compact in the
one cusped case.

Let $\xi$ be a totally irrational line starting from $e$. Fix a lift
$\widetilde{\xi}$ of $\xi$ to the universal cover $\widetilde{M}$ of
$M$, intersecting the preimages of the cut locus in successive points
$\widetilde{x_n}$.  Let $\widetilde{\sigma_n}$ be the lift of some
(simply connected by \ref{prop:cute_exist}) component of $\Sigma_0$
containing $\widetilde{x_n}$, and $\widetilde{r_n}$ the lift of $r_n$
starting from the same point at infinity $a$ as $\widetilde{\xi}$.

Since $r_n$ starts and ends in the cusp, then $\widetilde{r_n}$ ends
in $\gamma_n a$ for some $\gamma_n$ in the covering group of
$\widetilde{M}\ra M$.  By the definition of $r_n$, the geodesic ray
from $\widetilde{x_n}$ to $\gamma_n a$ projects to some minimizing ray
in $M$.

\myfigure{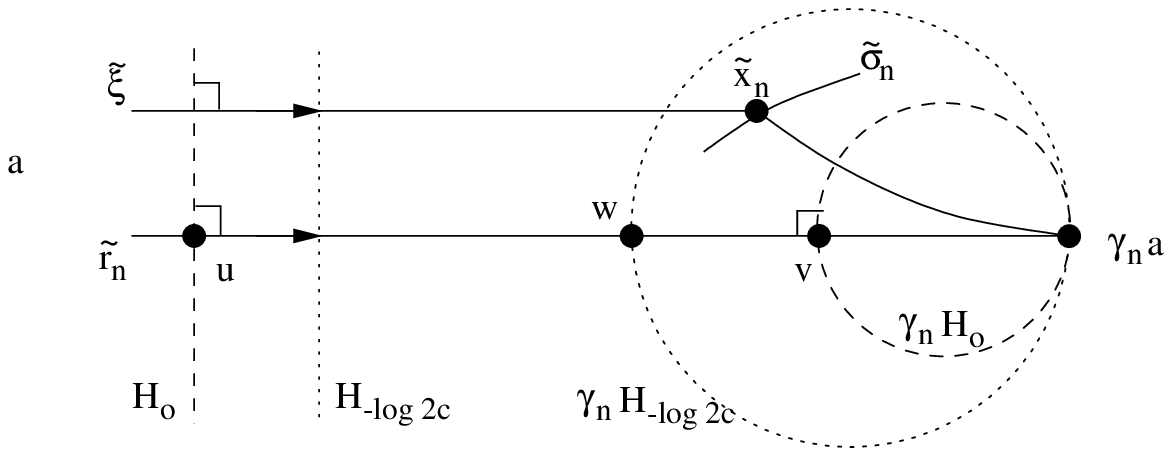}{4.2cm}{Convergence properties of the good 
approximating sequence}

Let $(H_t)_{t\in\RR}$ be the unique family of horospheres centered at 
$a$ such that $H_0$ maps onto $\beta_e^{-1}(0)$, the distance between 
$H_t$ and $H_{t'}$ is $|t-t'|$, and $H_t$ converges to $a$ as $t$ goes to 
$+\infty$. Since $-\log 2c = \min_{x\in\Sigma} \beta_e(x)$, each 
component of $\Sigma_0$ has a lift contained between the horospheres 
$H_0$ and $H_{-\log 2c}$. Hence $\widetilde{\sigma_n}$ is contained 
between $\gamma_n H_0$ and $\gamma_n H_{-\log 2c}$ (see Figure 
\arabic{fig}).

Let $u,w,v$ be the intersection points of $\widetilde{r_n}$ with
$H_0,\gamma_n H_{-\log 2c},\gamma_n H_0$, respectively. So that
$d_{\widetilde{M}}(u,v)=D(r_n)$ by definition of the depth, and
$d_{\widetilde{M}}(w,v)=|-\log 2c|=\log 2c$, since $\widetilde{r_n}$
is perpendicular to $\gamma_n H_0$ and $\gamma_n H_{-\log 2c}$. Since
$\widetilde{\xi}$ meets $\widetilde{\sigma_n}$, the horosphere
$\gamma_n H_{-\log 2c}$ which is centered at $\gamma_n a$ meets
$\xi$. Therefore, by the definition of the distance $d$ on $Lk(M,e)$,
with $\epsilon=+1$ if $w$ lies between $u$ and $v$, and $-1$
otherwise,
$$d(r_n,\xi) \leq \frac{1}{2} e^{-\epsilon d_{\widetilde{M}}(u,w)}
=\frac{1}{2} e^{-d_{\widetilde{M}}(u,v) + d_{\widetilde{M}}(w,v)}
=\frac{1}{2} e^{\log 2c} e^{-D(r_n)}.$$ This proves the result.
\eop{\ref{theo:good_indeed}}

\bcoro\label{coro:good_converges}
 Let $M$ be a non elementary geometrically finite pinched negatively
curved Riemannian manifold, with only one cusp $e$, having a cute cut
locus. The good approximating sequence $(r_n)_{n\in\NN}$ of any
totally irrational line $\xi$ starting from $e$ converges to $\xi$.  
\ecoro

\dem This follows from the previous theorem, and from the fact that
the depths $D(r_n)$ converge to $+\infty$. We prove this last claim by
absurd.  By Remark \ref{rem:non_numerote}, the sequence $(r_n)_{n\in
\NN}$ would have a constant subsequence $(r_{n_k})_{k\in \NN}$.  Let
$a$ be the starting point of $\widetilde{\xi}$ (with the notations of
the beginning of the previous proof), hence of $\widetilde{r_{n_k}}$.
Let $\Gamma_a \gamma_n \Gamma_a$ be the double coset associated to
$(r_n)_{n\in \NN}$, with $(\gamma_{n_k})_{k\in \NN}$ constant.

Since $\xi$ is irrational, and by strict convexity of the
horospheres, the sequence $\widetilde{\sigma_n}$ goes out of every
horosphere centered at $\gamma_{n_k}a$. This contradicts the fact that
$\widetilde{x_{n_k}}$ lies for each $k$ between the horospheres
$\gamma_{n_k}H_0$ and $\gamma_{n_k}H_{-\log 2c}$ (see Figure
\arabic{fig}).
\eop{\ref{coro:good_converges}}

\bigskip
Till the end of this section, we assume that $M$ is a non elementary
geometrically finite hyperbolic $3$-orbifold, uniformized as in
subsection \ref{sect:const_Curv_Case}. Again, this extends easily to
$\HH^{n}_{\RR}$.  We know a bit more on the depths of the good
approximants $r_n$, generalizing what Ford did in the case of the
Bianchi orbifold $\HH^3/PSL(2,\ZZ[i])$.

If $\xi$ is a totally irrational line starting from the cusp $e$, let
$\widetilde{\xi}$ be a lift of $\xi$ in $\HH^3$ starting from
$\infty$. It is a vertical geodesic in the upper halfspace $\HH^3$,
that cuts, while going downwards, a sequence
$(\widetilde{\sigma_n})_{n\in\NN}$ of cells of the preimage of
$\Sigma_0$ in $\HH^3$. The $n$-th cell
$\widetilde{\sigma_n}$ passed through by $\widetilde{\xi}$, is
contained in the boundary of the basins of $\gamma_n(\infty)$ and
$\gamma_{n+1}(\infty)$ for some $\gamma_n,\gamma_{n+1}$ in the
covering group $\Gamma$. By definition of the good approximation
sequence, the $n$-th good approximant $r_n$ is the rational
line associated to the double coset
$\Gamma_\infty\gamma_n\Gamma_\infty$.

\bprop
\label{lem:depth_increases}
If $(r_n)_{n\in\NN}$ is any good approximation sequence of a totally
irrational line $\xi$, then for all $n$ in $\NN$,
$$D(r_n) < D(r_{n+1}).$$
\eprop

\dem  With the notations above, by the definition of the summit
of $\widetilde{\sigma_n}$, there exist horospheres $H_n$ and $H_{n+1}$
centered respectively at $\gamma_n(\infty)$ and $\gamma_{n+1}(\infty)$
that are tangent at the summit $\widehat{\widetilde{\sigma_n}}$ of
$\widetilde{\sigma_n}$.  The Euclidean line through
$\widehat{\widetilde{\sigma_n}}$ perpendicular to
$\widetilde{\sigma_n}$ hence goes through the Euclidean center of the
hyperbolic plane $P$ containing $\widetilde{\sigma_n}$, and the
Euclidean centers of $H_n$ and $H_{n+1}$. Hence (see the Figure
below), the Euclidean radius of $H_{n+1}$, (which is contained in the
half-ball bounded by $P$ since $\widetilde{\xi}$ is first meeting
$B_{\gamma_n(\infty)}$ and then $B_{\gamma_{n+1}(\infty)}$), is
strictly smaller than the one of $H_n$.

\myfigure{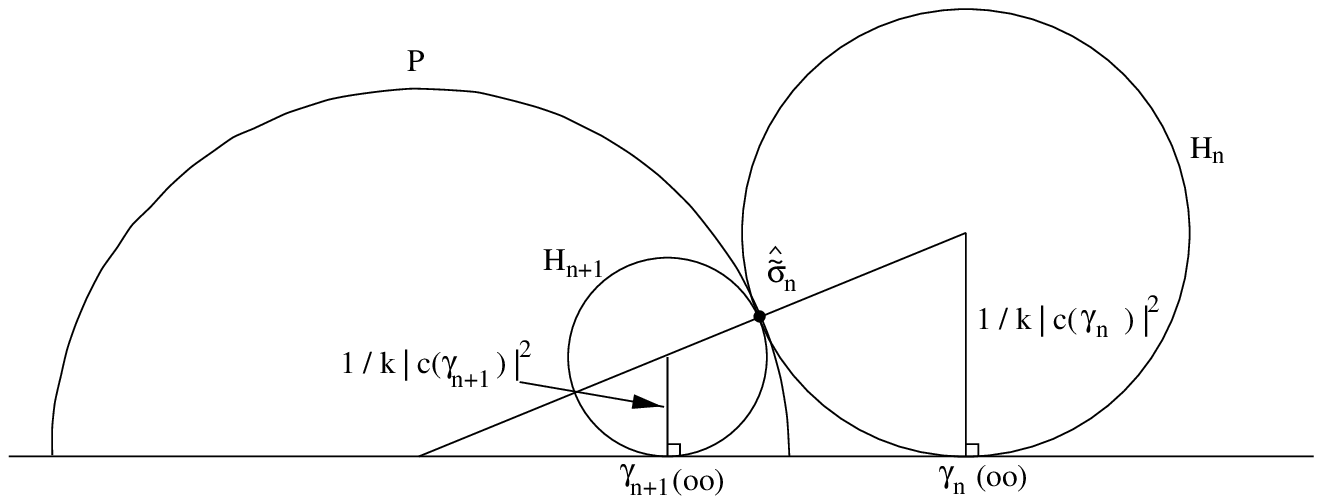}{4.5cm}{Depths of convergents are increasing}

Since the point $\widehat{\widetilde{\sigma_n}}$ is at the same
(hyperbolic) distance from the horospheres $H_{\gamma_n(\infty)}(1)$
and $H_{\gamma_{n+1}(\infty)}(1)$ centered at $\gamma_n(\infty)$ and
$\gamma_{n+1}(\infty)$ of Euclidean radii $\frac{1}{2|c(\gamma_n)|^2}$
and $\frac{1}{2|c(\gamma_{n+1})|^2}$, there is a constant $k>0$ such
that the Euclidean radii of $H_{n}$ and $H_{n+1}$ are respectively
$\frac{1}{2k|c(\gamma_n)|^2}$ and $\frac{1}{2k|c(\gamma_{n+1})|^2}$.
This proves the result, using Lemma \ref{lem:depth_computation}.
\eop{\ref{lem:depth_increases}}

\bigskip
Let ${\cal D}$ be the set of real numbers $e^{\frac{D(r)}{2}}$ for $r$
an integral line.  Since there are only finitely many integral lines,
the subset ${\cal D}$ of $[1,+\infty[$ is finite. It consists of
$\{1\}$ if there is only one integral line (counted without
multiplicity and without orientation), as in the case of the Bianchi
orbifolds $\HH^3/\,{\rm PSL}(2,{\cal O}_d)$ for $d=1,2,3,7,11$.

For $g,h$ elements of PSL$(2,\CC)$, define 
$$\Delta(g,h) =|\left|\begin{array}{cc} a(g)& a(h)\\ c(g) &
c(h)\end{array}\right| |.$$ This nonnegative real number does not
depend on the choosen lifts to SL$(2,\CC)$, nor on the representatives
of the left cosets of $g, h$ by the stabilizer of any horizontal
horosphere. In particular, by Lemma \ref{lem:depth_computation},
$$\Delta(1,h) = |c(h)|=e^{\frac{D(\Gamma_\infty h\Gamma_\infty)}{2}}$$
if $h\notin\Gamma_\infty$.

\bprop\label{prop:determinant_integral} 
If $r_n,r_{n+1}$ are consecutive good approximants of an irrational
line $\xi$, then $d(r_n,r_{n+1})\exp(\frac{D(r_n)+D(r_{n+1})}{2})$
belongs to the finite set ${\cal D}$.  In particular, if there is only
one integral line (counted without multiplicity and without
orientation), then
$$d(r_n,r_{n+1})=e^{ -(\frac{ D(r_n)+D(r_{n+1})}{2})}$$ 
\eprop

\dem We may assume that $\xi$ is totally irrational.  Choose
$\gamma_n,\gamma_{n+1}$ representatives of the double cosets
associated to $r_n,r_{n+1}$ such that the boundaries of the basins
$B_{\gamma_n(\infty)}$ and $B_{\gamma_{n+1}(\infty)}$ contain the
$n$-th cell of the preimage of $\Sigma_0$ passed through by
$\widetilde{\xi}$. By Lemma \ref{lem:depth_computation}, one has
$$\Delta(\gamma_n,\gamma_{n+1})= |\frac{a(\gamma_n)}{c(\gamma_n)}-
\frac{a(\gamma_{n+1})}{c(\gamma_{n+1})}|
|c(\gamma_n)||c(\gamma_{n+1})|=
d(r_n,r_{n+1})e^{\frac{D(r_n)+D(r_{n+1})}{2}}.$$ Since SL$_2(\CC)$
preserves the area in $\CC^2$,
$\Delta(\gamma_n,\gamma_{n+1})=\Delta(1,\gamma_n^{-1}\gamma_{n+1})$.
Since the closures of the basins $B_{\gamma_n(\infty)}$ and
$B_{\gamma_{n+1}(\infty)}$ meet in a cell of the preimage of
$\Sigma_0$, the rational line associated to the double coset of
$\gamma_n^{-1}\gamma_{n+1}$ is an integral line, and the result
follows.
\eop{\ref{prop:determinant_integral}}

\medskip
\rem When $M=\HH^2/{\rm PSL}_2(\ZZ)$, it is well known that
$\Delta(\gamma_n,\gamma_{n+1})$ is always $1$. This proposition gives
a geometric understanding of why this is not always the case, and a
geometric interpretation of the possible values, for instance for the
Poitou's best approximation sequence when $M=\HH^3/{\rm PSL}_2({\cal
O_{-d}})$ for $d$ large enough. See \cite{Poi} for the list of
possible values of $\left|\begin{array}{cc} p_n & p_{n+1} \\ q_n &
q_{n+1}\end{array}\right|$ when $d=19$ for example.

\section{Continued sequence of a totally  irrational line}
\label{sect:continued_fraction}

We keep the notations of the beginning of section
\ref{sect:rational_rays}. We assume that the cut locus of the cusp $e$
is cute. The aim of this section is to associate to a totally
irrational line $\xi$ starting from $e$ a sequence $(a_n)_{n\in \NN}$
in some countable alphabet which will determine it.

Let ${\cal R}$ be the finite set which consists of all the first
intersection points of integral lines starting from $e$ with
$L_e=\beta_e^{-1}(1)$. Let $\pi_1(L_e,{\cal R})$ denote the
set of homotopy classes relative to endpoints of paths in $L_e$ with
endpoints in ${\cal R}$. Note that unless there is only one point in
${\cal R}$, $\pi_1(L_e,{\cal R})$ is not a group, but it is a groupoid
for the composition of paths. There is a natural involution on
${\cal R}$. It associates to the first intersection point $\lambda$, of
an integral line $r$ with $L_e$ the second (and last) intersection
point, that we will denote by $\lambda^{-1}$. It is clear that
$\lambda^{-1}$ is also the first intersection point of the integral
line which is $r$ with the opposite orientation.  Since $L_e$ is
compact and ${\cal R}$ finite, the groupoid $\pi_1(L_e,{\cal R})$ is
countable.

Let $(r_n)_{n\in\NN}$ be the good approximation sequence of $\xi$,
with $(x_n)_{n\in\NN}$ the successive intersection points of $\xi$
with the cut locus and $\sigma_n$ the cell of $\Sigma_0$ containing
$x_n$. Recall that $\xi$ is oriented and totally irrational, hence
for each $n$, the tangent space to $\xi$ at $x_n$ is oriented and
transverse to the tangent subspace to $\sigma_n$ at $x_n$. For each
$n$, endow $\sigma_n$ with the transverse orientation given by the
oriented tangent space to $\xi$ at $x_n$.

\myfigure{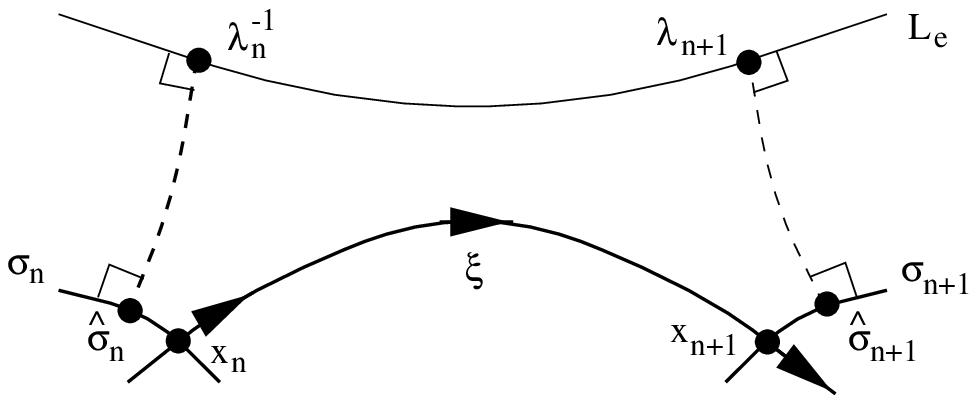}{3.5cm}{The continued  sequence
of an  irrational ray}

Note that $r_0$ is an integral line. Let $\lambda_0$ be its first
intersection point with $L_e$. Define $a_0\in \pi_1(L_e,{\cal R})$ to
be the class of the constant path at $\lambda_0$.  We define the
sequences $(a_n)_{n\in\NN},(\lambda_{n})_{n\in\NN}$ by
induction. Assume $a_n,\lambda_n$ are defined.  Consider the path $c$
starting at $\lambda_n^{-1}$, following the subpath of an integral
line from $\lambda_n^{-1}$ to the summit $\widehat{\sigma_n}$, then
any path $c_1$ in $\sigma_n$ from $\widehat{\sigma_n}$ to $x_n$, then
the subpath of $\xi$ from $x_n$ to $x_{n+1}$, then any path $c_2$ in
$\sigma_{n+1}$ from $x_{n+1}$ to $\widehat{\sigma_{n+1}}$, then the
subpath of an integral line starting from $\widehat{\sigma_{n+1}}$ on
the negative side of $\sigma_{n+1}$, that ends in the point
$\lambda_{n+1}\in{\cal R}$.  By pushing a little bit $c_1$ on the
positive side of $\sigma_n$, and $c_2$ on the negative side of
$\sigma_{n+1}$, one gets a path contained in the complement of the cut
locus in $M$. Recall that this complement canonically retracts onto
$L_e$. Let $a_{n+1}$ be the class of the path between $\lambda_n^{-1}$
and $\lambda_{n+1}$ obtained by retracting $c$ onto $L_e$.

\bdefi 
The sequence $(a_n)_{n\in\NN}$ is called the continued sequence
of the totally irrational ray $\xi$.  
\edefi

\btheo
\label{theo:continued_fraction_determines_ray}
Let $M$ be a non elementary geometrically finite pinched negatively
curved Riemannian manifold, with only one cusp $e$, having a cute cut
locus. A totally irrational line is uniquely determined by its continued 
sequence.
\etheo

\dem Let $(a_n)_{n\in\NN}$ be the continued sequence of the totally
irrational ray $\xi$. Let $\lambda_{n} ^{-1}$ be the initial point of
the path $a_n$. Consider the path $v$ in $M$, depending only on the
continued sequence, which is obtained by following the minimizing
geodesic ray starting from $e$ arriving perpendicularly on
$\lambda_0$, then the subpath of the integral line from $\lambda_0$ to
$\lambda_0^{-1}$, then $a_1$, then the subpath of the integral line
from $\lambda_1$ to $\lambda_1^{-1}$, then $a_2$, etc. It is clear by
the construction that $v$ is homotopic to $\xi$, by an homotopy which
is proper on each negative subray (though not necessarily globally
proper).  Since $\Sigma$ is compact, the sequence of points
$(\lambda_{n})_{n\in\NN}$ on $v$ stays at uniformely bounded distance
from $\xi$. If $\xi'$ is another totally irrational line having the
same continued sequence, lift to the universal cover the homotopy
between $v$ and $\xi$ and the one between $v$ and $\xi'$, so that they
coincide on some lift of $v$. Then the lifts of $\xi,\xi'$ are two
geodesic lines in $\widetilde{M}$ that have the same origin at
infinity, and a sequence of point converging to their endpoint at
infinity that stay at uniformely bounded distance.  Hence the geodesic
lines have the same endpoint at infinity, therefore they coincide.
By projecting to $M$, one gets that $\xi,\xi'$ are equal. 
\eop{\ref{theo:continued_fraction_determines_ray}}

\bigskip
Till the end of this section, we assume that $M$ is a non elementary
geometrically finite hyperbolic $3$-orbifold, uniformized as in
subsection \ref{sect:const_Curv_Case}.  The orbifold universal cover
of $L_e$ is the horizontal horosphere in $\HH^3$ which is mapped onto
$L_e$ by the choice of the orbifold universal cover $\HH^3\ra M$. We
will identify that horosphere with $\CC$ by vertical projection. Let
$\widetilde{{\cal R}}$ be the subset of $\CC$ corresponding to the
lift of ${\cal R}$, which is a discrete subset of $\CC$. It is
invariant by the group $\Gamma_\infty$ (and it is reduced to one orbit
if there is only one integral line (counted without multiplicity and
orientation)).

Assume for simplicity that the group $\Gamma_\infty$ is a covering 
group of $L_e$, hence a group of translations of $\CC$. Since any 
path between two points in $\widetilde{{\cal R}}$
is homotopic relative to endpoints to the segment between the 
endpoints, one can naturally identify $\pi_1(L_e,{\cal R})$
with the set $\widetilde{{\cal R}}-\widetilde{{\cal R}}$ of differences 
of two elements of $\widetilde{{\cal R}}$. Indeed,  any path
between two points of ${\cal R}$ has a unique lift, once a preimage of
the starting point is choosen, between two points in $\CC$, and we associate
to the path the difference of the endpoints of the lift. This does not 
depend on the  choosen preimage of the 
starting point, since any two of them differ by an element of 
$\Gamma_\infty$, which acts by translation.

Let $\xi$ be a totally irrational line, and $(r_n)_{n\in\NN}$ be the
good approximating sequence for $\xi$. Let $\widetilde{\xi}$ be a lift
of $\xi$ starting from $\infty$, and $\widetilde{r_n}$ the lift of
$r_n$ obtained by lifting the homotopy between $r_n$ and $\xi_n$ (the
subpath of $\xi$ up to the $n$-th intersection point $x_n$ with the
cut locus, followed by the minimizing geodesic ray from $x_n$ to $e$
on the opposite side). Let $\Gamma_\infty\gamma_n\Gamma_\infty$ be the
double coset associated to $r_n$, with $\gamma_n$ a representative so
that the endpoint of $\widetilde{r_n}$ is $z_n=\gamma_n(\infty)$ (this
determines the left coset $\gamma_n\Gamma_\infty$). Define by
convention $\gamma_{-1}=id$, so that $z_{-1}=\infty$.

The next result explains how the continued fraction can be
explicitely computed in terms of the good approximation sequence.

\bprop
\label
{prop:continued_sequence_constcurv}
Under the above identification of $\pi_1(L_e,{\cal R})$
with $\widetilde{{\cal R}}-\widetilde{{\cal R}}$, one has
 $a_0=0$ and for $n\geq 0$,
$$a_{n+1} = \gamma_n^{-1}(z_{n+1})-\gamma_n^{-1}(z_{n-1}).$$ 
\eprop

\noindent Note that the right hand side does not depend on the left
coset of $\gamma_n$.

\medskip
\dem Let $\widetilde{x_n},\widetilde{x_{n+1}}$ be consecutive
intersection points of $\widetilde{\xi}$ with the preimage of the cut
locus. Then $\widetilde{\xi}$ passes at $\widetilde{x_n}$ from the
basin of $\gamma_{n-1}(\infty)$ to the basin of $\gamma_{n}(\infty)$,
stays inside the basin of $\gamma_{n}(\infty)$ between
$\widetilde{x_n}$ and $\widetilde{x_{n+1}}$, then passes at
$\widetilde{x_{n+1}}$ into the basin of
$\gamma_{n+1}(\infty)$. Consider the action of $\gamma_n^{-1}$ on
$\HH^3$. It maps $z_n$ to $\infty$ and $z_{n\pm 1}$ to
$\gamma_n^{-1}(z_{n\pm 1})$. It preserves the set of lifts of any
integral line. Since the geodesic lines starting from $\infty$ are
vertical (half-)lines, the summits of the cells
$\gamma_n^{-1}(\widetilde{\sigma_n})$ and
$\gamma_n^{-1}(\widetilde{\sigma_{n+1}})$ projects vertically to
$\gamma_n^{-1}(z_{n- 1})$ and $\gamma_n^{-1}(z_{n+1})$ respectively.
The result follows.  \eop{\ref{prop:continued_sequence_constcurv}}

\bcoro\label{coro:cont_seq_constcurv} Assume that there is only one
integral line (counted without multiplicity). Then
$$|a_{n+1}| = \Delta(\gamma_{n+1},\gamma_{n-1}).$$ 
\ecoro

\dem 
$$\Delta(\gamma_{n+1},\gamma_{n-1})=
\Delta(\gamma_n^{-1}\gamma_{n+1},\gamma_n^{-1}\gamma_{n-1})
=|\gamma_n^{-1}\gamma_{n+1}(\infty)-\gamma_n^{-1}\gamma_{n-1}(\infty)|
= |a_{n+1}|.$$
 \eop{\ref{coro:cont_seq_constcurv}}

\medskip
The next result proves that the good approximation sequence can be
recovered from the continued sequence. To get $r_n$, one only has to
compute the value of the endpoint $z_n\in\CC$ of the lift
$\widetilde{r_n}$. We prove that $z_n$ can be expressed in an explicit
summation formula in terms of the $a_i$'s with $i\leq n$ and the
denominators of the $z_i$'s with $i<n$.

\btheo\label
{theo:new_formula}
If $\xi$ is a totally irrational line, then with the above notations and 
$q_i=c(\gamma_i)$ for $i\geq -1$, one has
$$z_{n}=  z_0+\sum_{k=1}^{n} \frac{1}{\sum_{i=0}^{k}
(-1)^iq_{i-1}^2 a_i}.$$
\etheo

\dem
By Proposition \ref{prop:continued_sequence_constcurv}, one has
$$z_{n+1}=\gamma_n(a_{n+1} + 
\gamma_n^{-1}z_{n-1}).$$
Writing for simplicity $\gamma_n=\left(\begin{array}{cc} a & b \\ c & 
d\end{array}\right)$, one gets, since $ad-bc =1$, for $n\geq1$
$$z_{n+1}=\frac{a\left(a_{n+1} +
\frac{d z_{n-1}-b}{-c z_{n-1}+a}\right)+b}{c
\left(a_{n+1} + 
\frac{d z_{n-1}-b}{-c z_{n-1}+a}\right)+d}$$
$$=\frac{z_{n-1} + a^2 a_{n+1} - ac a_{n+1}z_{n-1}}
{aca_{n+1} -c^2a_{n+1}z_{n-1} +1}$$
$$= \frac{z_{n-1} + 
z_{n}c^2a_{n+1}(z_{n}-z_{n-1})}{1+
c^2a_{n+1}(z_{n}-z_{n-1})},$$
since $z_{n}=\gamma_n(\infty)=\frac{a}{c}$. Hence denoting $q_n=c(\gamma_n)$,
which depends (up to sign) only on $r_{n}$, one gets 
$$
z_{n+1}= \frac{z_{n-1} + 
z_{n}q_n^2a_{n+1}(z_{n}-z_{n-1})}{1+
q_n^2a_{n+1}(z_{n}-z_{n-1})}.
$$
Hence (upon adding and subtracting $z_n$ in the numerator) we have 
$$z_{n+1}-z_{n}= -\frac{z_{n}-z_{n-1}}{1+
q_n^2a_{n+1}(z_{n}-z_{n-1})}.$$ Let $x_{p}=\frac{1}{z_{p}-z_{p-1}}$
for $p\geq 1$ and $x_0=0$.  Then for $n\geq 0$
$$x_{n+1}= -(x_n +q_n^2a_{n+1}).$$ Since by convention $q_{-1}=0$, one
gets $x_n=\sum_{k=0}^{n} (-1)^k q_{k-1}^2a_k$.  The equality of the 
theorem now
follows from the fact that $z_{n+1}-z_{n}=\frac{1}{x_{n+1}}$.
\eop{\ref{theo:new_formula}}

\bigskip
In the case of $\HH^2/\PSLZ$, our continued sequence $(a_n)_{n\in\NN}$
of $\xi$ slightly differs from the classical continued fraction expansion
$$\frac{1}{b_1+\frac{1}{b_2+\frac{1}{...}}}$$ of the endpoint on the
real axis of the lift of $\xi$ starting at $\infty$ and ending in
$]-1,1[$: one has $|a_n-b_n|\leq 1$ (and the difference may be a
non recursive function).  This is due to working with the cut locus
rather than with its dual cell decomposition, but we will come back to
that point in another paper.

\bigskip
\noindent {\small\begin{tabular}{ll}
\begin{tabular}{l} 
Caltech \\
Department of Mathematics 
\\ Pasadena CA 91125, USA \\
{\it e-mail: saar@cco.caltech.edu}
\end{tabular}
&
\begin{tabular}{l}
{\it Current address:}\\ Department of Mathematics\\ 
University of Ben Gurion \\ BEER-SHEVA, Israel\\
{\it e-mail: saarh@cs.math.bgu.ac.il}
\end{tabular}
\end{tabular}
\\
 \mbox{}
\\
 \mbox{}
\\
\hspace*{0.1cm} 
\begin{tabular}{l} Laboratoire de Math\'ematiques UMR 8628 CNRS\\
Equipe de Topologie et Dynamique (B\^at. 425)\\
Universit\'e Paris-Sud \\
91405 ORSAY Cedex, FRANCE.\\
{\it e-mail: Frederic.Paulin@math.u-psud.fr}
\end{tabular}
}


\begin{thebibliography}{99}{\small


\bibitem[Ah]{Ah}
L.V.~Ahlfors, {\it M\"obius transformations and Clifford numbers},
Differential Geometry and Complex Analysis-in memory of H.E.~Rauch (I.~
Chavel and H. Farkas, eds.), Springer-\-Verlag, New York, 1985.

\bibitem[Ano]{Ano}
D.V.~Anosov, {\it Geodesic flows on closed Riemann manifolds with negative 
curvature}, Proc. Steklov Inst. Math. {\bf 90} (1967), Amer. Math. Soc., 
1969.

\bibitem[BGS]{BGS}
 W.~Ballmann, M.~Gromov, V.~Schroeder,  {\it Manifolds of nonpositive 
curvature}, Progress in Math. {\bf 61}, Birkh\"auser 1985.

\bibitem[Bea]{Bea}
A.~Beardon, {\it The geometry of discrete groups}, Springer-Verlag, 1983.

\bibitem[Bow]{Bow} B.~Bowditch, {\it Geometrical finiteness with
variable negative curvature}, Duke Math. J. {\bf 77} (1995), 229-274.

\bibitem[BK]{BK}
P.~Buser, H.~Karcher, {\it Gromov's almost flat manifolds}, 
Ast\'erisque {\bf 81}, Soc.~Math.~France, 1981.

\bibitem[Dan]{Dan}
S.J.~Dani, {\it Bounded orbits  of flows on homogeneous spaces},
Comment. Math. Helv. {\bf 61} (1986), no. 4, 636--660. 

\bibitem[DP]{DP} 
R.~Descombes, G.~Poitou, {\it Sur l'approximation dans $\RR(i\sqrt{11})$}, 
Compt. Rend. Acad. Scien. Paris {\bf 231} (1950) 264-264.

\bibitem[EP]{EP} 
D.B.A.~Epstein, R.~Penner, {\it Euclidean decomposition of non-compact 
hyperbolic manifolds}, J. Diff. Geom. {\bf 27} (1988),
67-80.

\bibitem[FLP]{FLP} 
A.~Fathi, F.~Laudenbach, V.~Poenaru, {\it Travaux  
de Thurston sur les surfaces}, Ast\'erisque {\bf 66-67}, 
Soc. Math. France 1979.

\bibitem[Fel]{Fel}
E.A.~Feldman,  {\it The geometry of immersions I},  
Trans. Amer. Math. Soc. {\bf 120} (1965) 185-224.

\bibitem[For1]{For}
L.~Ford, {\it Rational approximations to irrational complex numbers}, 
Trans. Amer. Math. Soc. {\bf 99} (1918), 1-42.

\bibitem[For2]{For2}
L.~Ford, {\it On the closeness of approach of complex rational fraction
to a complex irrational number}, 
Trans. Amer. Math. Soc. {\bf 27} (1925), 146-154.

\bibitem[GH]{GH}
E.~Ghys, P.~de la Harpe, eds. {\it Sur les groupes hyperboliques d'apr\`es 
Mikhael Gromov}, Prog. in Math. {\bf 83}, Birkh\"auser 1990.

\bibitem[Haa]{Haa}
A.~Haas, {\it Diophantine approximation on hyperbolic surfaces}, 
Acta Math. {\bf 156} (1986) 33-82.

\bibitem[HS]{HS}
A.~Haas-C.~Series, {\it The Hurwitz constant and diophantine 
approximation on Hecke groups}, 
J. Lond. Math. Soc. {\bf 34} (1986) 219-234.

\bibitem[Ham]{Ham}
U.~Hamenst\"adt, {\it A new description of the Bowen-Margulis measure},
Erg.~Theo.~Dyn.~Sys. {\bf 9} (1989)  455--464. 

\bibitem[Hat]{Hat} 
A.~Hatcher, {\it Hyperbolic structures of arithmetic type on some 
link complements}, 
J. Lond. Math. Soc. {\bf 27} (1983) 345-355.



\bibitem[HP]{HP}
S.~Hersonsky, F.~Paulin, {\it  On the rigidity of discrete isometry 
groups of negatively curved spaces}, 
Comm. Math. Helv. {\bf 72} (1997) 349-388.

\bibitem[HV]{HV}
R.~Hill-S.L.~Velani, {\it The Jarn\'{\i}k-Besicovitch theorem for
geometrically  finite Kleinian groups},
Proc. London Math. Soc. {\bf 3} (1997) 524-551.
 
\bibitem[Hof]{Hof36}
N.~Hofreiter, {\it Diophantische  Approximationen in imagin\"aren 
quadratischen Zahlk\"orpern}, 
Monatsh. Math. Phys {\bf 45} (1936), 175-190.


\bibitem[LS]{LS}
J.~Lehner, M.~Sheingorn, {\it Simple closed geodesics on $H_{+}/\Gamma (3)$ arise from the Markov spectrum}, 
Bull. Amer. Math. Soc. {\bf 11} (1984) 359-362.

\bibitem[Pat]{Pat}
S.J.~Patterson, {\it Diophantine approximation in Fuchsian groups},
Philos. Trans. Roy. Soc. London Ser. A {\bf 282}, 241-273.

\bibitem[Per1]{Per31}
O.~Perron, {\it \"Uber einen Approximationssatz von Hurwitz und \"uber die 
Approximation einer komplexen Zahj durch Zahlen des K\"orpers der dritten 
Einheitswurzeln}, 
S.-B.~Bayer Akad. Wiss. (1931), 129-154.

\bibitem[Per2]{Per33}
O.~Perron, {\it Diophantische  Approximationen in imagin\"aren 
quadratischen Zahlk\"orpern, insbesondere im K\"orper $\RR(i\sqrt{2})$}, 
Math. Z. {\bf 37} (1933), 749-767.

\bibitem[Poi]{Poi}
G.~Poitou, {\it Sur l'approximation des nombres complexes par les nombres 
des corps imaginaires quadratiques d\'enu\'es d'id\'eaux principaux}, 
Ann. Scien. Ec. Norm. Sup. {\bf 70} (1953), 199-265.

\bibitem[Rat]{Rat}
J.G.~Ratcliffe, {\it Foundations of Hyperbolic Manifolds},
Springer-Verlag, New-York, 1994.

\bibitem[Schmi]{Sch} A.L. Schmidt, {\it Diophantine approximations of
complex numbers}, Acta. Math. {\bf 134} (1975), 1-84.

\bibitem[Schmu]{Schmu}
P.~Schmutz-Schaller, {\it The modular torus has maximal length spectrum},  
GAFA {\bf  6} (1996)  1057-1073.

\bibitem[Ser1]{Ser1} C.~Series, {\it On coding geodesics with
continued fractions}, L'Enseign. Math., {\bf 29} (1980) 67-76.

\bibitem[Ser2]{Ser2} C.~Series, {\it Symbolic dynamics for geodesic
flows}, Acta Math., {\bf 146} (1981) 103-128.

\bibitem[Ser3]{Ser3}
C.~Series,  {\it The modular surface and continued fractions}, 
J.~Lond. Math. Soc. {\bf 31} (1985) 69-80.

\bibitem[Sug]{Sug} 
K.~Sugahara, {\it On the cut locus and the topology of manifolds},
J.~Math.~Kyoto Univ. {\bf 14} (1974), 391-411.

\bibitem[Sul]{Sul} 
D.~Sullivan, {\it Disjoint spheres, approximation
by imaginary quadratic numbers, and the logarithm law for geodesics},
Acta Math., {\bf 149} (1982), 215-237.

\bibitem[Swa]{Swa} 
R.~Swan, {\it Generators and relations for certain special linear groups}, 
Adv. Math. {\bf 6} (1971) 1-77.

\bibitem[Tho]{Tho}
R.~Thom, {\it La stabilit\'e topologique des applications polynomiales}, 
L'Ens. Math.   {\bf 8} (1962), 24-33.

\bibitem[Tro1]{Tro1} 
D.~Trotman, {\it Interpr\'etation topologique des conditions de Whitney}, 
Ast\'erisque {\bf 59-60} Soc. Math. France 1978, 233-248.

\bibitem[Tro2]{Tro} 
D.~Trotman, {\it Equisingularit\'e et conditions de Whitney}, 
Th\`ese, Universit\'e d'Orsay (1980).

\bibitem[Vul1]{Vul1}
L.~Vulakh, 
{\it Diophantine approximation on Bianchi groups},
J. Number Theo. {\bf 54} (1995) 73-80.
  
\bibitem[Vul2]{Vul2}
L.~Vulakh, 
{\it Diophantine approximation in $\RR^n$},
Trans. Amer. Math. Soc. {\bf 347} (1995) 573-585.

\bibitem[Whi]{Whi}
H.~Whitney, {\it Local properties of analytic varieties}, 
Differential and Combinatorial Topology, S.~Cairns ed.  (1965) 205-244. 

 }
\end{thebibliography}
\end{document}